\documentclass{article}
\usepackage{amsfonts}
\usepackage{latexsym}
\usepackage{amssymb}
\usepackage{amsmath}
\usepackage{mathrsfs}
\usepackage[mathscr]{euscript}
\usepackage{mathrsfs}
\usepackage{marvosym}
\usepackage{amsmath}
\usepackage{amssymb}
\usepackage{fancyhdr}
\usepackage{cite}
\usepackage{graphicx}
\usepackage{amsmath}
\usepackage[numbers,sort&compress]{natbib}
\usepackage{hyperref}
\usepackage{color}
\usepackage{comment}

\newtheorem{lemma}{Lemma}[section]
\newtheorem{corollary}{Corollary}[section]
\newtheorem{proposition}{Proposition}[section]
\newtheorem{theorem}{Theorem}[section]
\newtheorem{remark}{Remark}[section]
\newtheorem{definition}{Definition}[section]
\newtheorem{example}{Example}[section]
 
\def\cvd{\hfill $\Box$}
\def\R{\mathop{{\rm I}\mskip-4.0mu{\rm R}}\nolimits}

\title{{{\bf\large {\color{black}{New Separation Theorems between Convex Sets in Groups and Applications}}}}\thanks{This
work was supported by the National Natural Science Foundation of China
(11871059), the Innovation Team Funds of China West Normal University (KCXTD2023-3) and the Meritocracy Research Funds of China West Normal University (17YC379).}}

\author {\small Jun Li\thanks{Laboratory of Optimization Theory and Applications, School of Mathematical Sciences, China West Normal University, 637009 Nanchong, Sichuan, P. R. China (Corresponding author, e-mail: junli@cwnu.edu.cn)}\, and\, Giandomenico Mastroeni \thanks{Department of Computer Science, University of  Pisa, Largo B. Pontecorvo 3, 56127 Pisa,  Italy (gmastroeni@di.unipi.it)}
}

\begin{document}
\date{}
\maketitle{}

{\bf Abstract.} The aim of this paper is to prove several new separation theorems between convex sets in groups under the assumption that
at least one of them has nonempty pseudo relative core, whereas their quasi cores are empty,  in which we  apply an extension theorem for  additive functions  but not the Hahn-Banach theorem for groups. This answers affirmatively the open question proposed by Li and Mastroeni recently. As an application, we provide Lagrangian type necessary optimality conditions  for a class of optimization problems with cone constraints in groups and semigroups.



\def\cvdeq{\hfill \sqcap \hskip-5.5pt \sqcup}

\def\vs{\vskip.3truecm}

\def\inte{{\rm \ int\ }}
\def\cl{{\rm \ cl\ }}

\section{Introduction}
\setcounter{equation}{0}
Many cornerstone results in functional analysis, such as the Hahn-Banach theorem, are naturally expressed and proved via convexity.
Convex analysis is of fundamental importance in the analysis of optimization, variational inequalities  and monotone operator theory \cite{BP86, BC17, BL06, BV10, HL93, Roc70, rw98, sw70, Zal02}.
As is well known, separation theorem has been proved to be one of the most fertile one in convexity theory and their applications.

Recently, Borwein  and Giladi \cite{BG18} proposed a  convexity notion for groups and monoids in the field of non-divisible groups, which collapses to
the classical one when applied to a vector space. They proved that several  results of classic convex analysis  also hold for functions defined on
groups and semigroups and provided suitable examples of applications.
Moreover, they further proved separation  and minimax theorems  \cite{BG16} in the particular setting of locally convex topological groups and monoids.

In \cite{LM20}, Li and Mastroeni proved the equivalence of the definitions of convex sets and  functions  introduced by Adivar  and Fang, Borwein and Giladi in finite-dimensional  integer spaces.
   They further introduced the definition of  convex cones, affine sets and  analysed the separation between convex sets.
Recently, Li and Mastroeni \cite{LM24} investigated the main properties of convex sets, convex cones, and convex functions in groups and semigroups.
In particular, they introduced the concept of the quasi core of a set to prove  separation theorems between two convex subsets of a commutative group under the assumption of nonemptyness of the  quasi core of at least one of them. They also studied optimization problems defined on groups and semigroups by using the obtained results \cite{LM23, LM24}.
Very recently, Li and Mastroeni \cite{LM25} introduced the  pseudo relative core of a set in groups and investigated many interesting,  elegant and  useful properties of this new notion.
As applications, they extended the classical Radon's lemma and Helly's theorem from finite-dimensional Euclidean spaces to $\mathbb{Z}^d$ and groups and they analyse optimality conditions for a class of optimization problems with cone constraints in groups. They also proposed a question is: Can we show the (proper) separation between convex sets in groups under the assumption that
at least one of them has nonempty pseudo relative core, whereas their quasi cores are empty, or without any other additional assumptions?

In this paper, we prove several new separation theorems between convex sets in groups by using the extension theorem for additive functions in groups.  As an application, we provide Lagrangian type necessary optimality conditions  for a class of optimization problems with cone constraints in groups and semigroups.
Our main contributions are the following:
\begin{itemize}
\item  We provide a refinement of generalised affine functions on commutative groups. Specifically, we prove that a generalised affine function on commutative groups can not attain values $+\infty$ and $-\infty$ at the same
time. Moreover, we characterize additive functions on commutative groups and present some relationships between  additive, subadditive and generalised affine functions on commutative groups.

\item  We prove several new proper separation theorems between convex sets in commutative groups under the assumption that
at least one of them has nonempty pseudo relative core by using the extension theorem for additive functions in groups.
Compared with {\rm \cite[Theorem 5.3]{LM24}}, we remove the assumption that at least one of them has nonempty quasi core, where the Hahn-Banach theorem for groups {\rm \cite[Theorem 7]{BG18}} was employed. This answers affirmatively the open question proposed by Li and Mastroeni in {\rm\cite{LM25}}.
	
\item We provide Lagrangian type necessary optimality conditions for a class of optimization problems with cone constraints in groups and semigroups  under the assumption that the involved functions are convex-like. Compared with {\rm \cite[Theorems 6.2 and 6.3]{LM25}}, we do not assume ${\rm qcore (cone}\,{\mathcal{E}}_{\bar x})\ne\emptyset$.
\end{itemize}

The paper is organized as follows. In section 2, we  recall the preliminary concepts and notations  used throughout the paper.
In section 3, we investigate convex sets and cones in groups and semigroups, while, in section 4, we recall the main properties of  additive, subadditive,generalised affine, convex and convex-like functions  in groups.
In section 5, we present separation theorems between convex sets in groups by using the extension theorem  for additive functions in groups.
As an application, we provide  necessary optimality conditions  for a class of optimization problems with cone constraints in groups and semigroups in section 6.

\section{Basic concepts and notations}
%\indent
 Let $\mathbb{R}^m$ be the $m$-dimensional Euclidean space, where $m\in \mathbb{N}$, the positive integers.  Let $\mathbb{R}^m_+:=\{x:=(x_1,\cdots, x_m)\in \mathbb{R}^m: {\rm each }\,\, x_i\geq 0\}$,
%$\mathbb{R}^m_{++}:=\{x:=(x_1,\cdots, x_m)\in \mathbb{R}^m: {\rm each }\,\, x_i> 0\}$,$\mathbb{Q}^m:=\{x:=(x_1,\cdots, x_m) \in \mathbb{R}^m: {\rm each }\,\, x_i\in\mathbb{Q}\}$,
$\mathbb{Z}^m:=\{x:=(x_1,\cdots, x_m) \in \mathbb{R}^m: {\rm each }\,\, x_i\in\mathbb{Z}\}$, and $\mathbb{Z}^m_+:=\{x:=(x_1,\cdots, x_m) \in \mathbb{Z}^m: {\rm each }\,\, x_i\geq 0\}$.
We denote   by $\mathbb{Q}_+$  the non-negative rational numbers and by $\|\cdot\|$ the  Euclidean 2-norm in $\mathbb{R}^m$. All vectors  are considered as row vectors.

 We recall in this section the main concepts and notations for groups and semigroups.

\begin{definition}\label{Def:2.1} {\rm  (Semigroup and group) }  {\it  A set $X$ endowed with a binary operation $\cdot$ is called a
semigroup, if  the associative law $(a\cdot b)\cdot c = a\cdot(b\cdot c)$ is satisfied for all $a, b, c\in X$.
$X$ is called commutative, if $a\cdot b=b\cdot a$ for all $a, b \in X$. $e\in X$ is called a
unit element in $X$, if it satisfies $e\cdot x=x\cdot e=x$ for all $x\in X$. A set $X$ endowed with a binary operation $\cdot$ is called a group, if it is a semigroup with the unique unit element $e\in X$, and the unit inverse of each $x\in X$, say $x^*$, i.e., $x^*\cdot x=x\cdot x^*=e$.
}
\end{definition}

\begin{definition}\label{Def:2.2} {\rm  (Monoid) }  {\it  A set $X$ endowed with a binary operation $\cdot$ is called a
monoid, if $(X,\cdot)$ is a commutative semigroup with unit element.
}
\end{definition}

\begin{definition}\label{Def:2.3} {\rm  (Semiring and ring) }  {\it  A set $X$ endowed with two binary operations $+$ and $\cdot$ is called
a semiring (resp., ring), if $(X,+)$ is a commutative semigroup (resp., group) with unit element $0$, and $(X, \cdot)$ is
a semigroup (not necessarily commutative) with unit element $1$. The two operations
are connected by the distributive laws. We additionally assume for a semigroup $X$ that $r0=0r=0$  for all $r \in X$.
}
\end{definition}

\begin{definition}\label{Def:2.4} {\rm  (Semimodule and module)}  {\it A set $X$ endowed with a binary operation $+$ is called
 a left semimodule (resp., module) over a semiring (resp., ring) $R$, if $(X,+)$ is a commutative semigroup (resp., group) with
unit element $0$, and there is a map $R\times X\to X$, $(r, x) \mapsto rx$ such that for all
$r, r^1, r^2\in R$ and $x, x^1, x^2\in X$,
\begin{itemize}
\item[\rm(i)]  $(r^1r^2)x = r^1(r^2x)$.

\item[\rm(ii)] $r (x^1 + x^2) = rx^1 + rx^2$.

\item[\rm(iii)] $(r^1 + r^2)x = r^1x + r^2x$.

\item[\rm(iv)] $1x = x$.

\item[\rm(v)] $r0 = 0 = 0x$.
\end{itemize}
}
\end{definition}

We observe that a monoid is a semimodule over the non-negative
integers $\mathbb{Z}_+$ and  a commutative group is a module over the ring $\mathbb{Z}$  and also a monoid (i.e., a semimodule over the non-negative integers $\mathbb{Z}_+$).

%\begin{definition}\label{Def:2.5}{\rm  (Topological group and monoid) }  {\it A group endowed with a topology is said to be a
%topological group if the group operations are continuous. That is, the function $(x,y)\mapsto x-y$   is continuous. A monoid which is also a topological space is said to
%be a topological monoid if the addition operation is continuous.
%}
%\end{definition}

\begin{definition}\label{Def:2.6}{\rm  (Semidivisible, divisible and uniquely  divisible semimodule (resp., module)) }  {\it Let $X$ be a semimodule (resp. module) over a semiring (resp., ring) $R$.
\begin{itemize}
\item[\rm(i)] $X$ is called divisible, if for every $y\in X$ and $r\in R\setminus\{0\}$, there is $x \in X$ such that $rx=y$.

\item[\rm(ii)] $X$ is called uniquely  divisible, if for every $y\in X$ and $r\in R\setminus\{0\}$, there is a unique $x\in X$ such that $rx=y$.

\item[\rm(iii)] Let $X$ be a monoid. $X$ is called  $p$-semidivisible, if there is $p\in \mathbb{N}$ prime such that $pX=X$, and $X$ is called semidivisible, if it  is $p$-semidivisible for some prime $p$.
\end{itemize}
}
\end{definition}

\begin{definition}\label{Def:2.7-}{\rm  ($n$-generator and cyclic group) }  {\it Let $A$ be a nonempty subset of a group $X$ and $n\in \mathbb{N}$. The subgroup generated by $A$, denoted by $\langle A \rangle$, is the intersection of all subgroups of $X$ which contain $A$. A group is said to be an $n$-generator group, if it can be generated by some $n$-subset $\{x^1,\cdots, x^n\}$. A group is finitely generated, if it is $n$-generator  for some $n$.
A $1$-generator group $\langle x \rangle:=\langle\{x\}\rangle$ is said to be cyclic.
}
\end{definition}

If $\{X^\lambda:\lambda\in \Lambda\}$ is the set of subgroups of a group $X$, the subgroup generated by the $X^\lambda$' s is defined to be $\langle\bigcup_{\lambda\in \Lambda}X^\lambda\rangle$. This will be written
$\langle X^\lambda:\lambda\in \Lambda\rangle$ or in case $\Lambda:=\{\lambda^1,\cdots, \lambda^n\}$, a finite set, $\langle X^{\lambda^1},\cdots, X^{\lambda^n}\rangle$, where $n\in \mathbb{N}$.

\begin{definition}\label{Def:2.7}{\rm  (Finite and infinite order) }  {\it Let $X$ be a group. A element $x\in X$ is said to be have finite order $n\in \mathbb{N}$, if the cyclic subgroup $\langle x\rangle$ of $X$ has order $n$.
A element $x\in X$ is said to be have infinite order, if the cyclic subgroup $\langle x\rangle$ of $X$ has infinite order.
}
\end{definition}

\begin{lemma} \label{LM:2.2}  {\rm \cite[1.3.8]{Rob96}} {\it Let $x$ be a element of a group $X$. $x$ has infinite order if and only if $mx\neq nx$ for $m\neq n$, where $m, n\in \mathbb{N}$, or equivalently,  $kx^0=0$ implies that $k=0$, where $k\in \mathbb{Z}$.
}
\end{lemma}

\begin{definition}\label{Def:CPO} {\rm  (Compatible partial order on a semimodule) }  {\it Given a semimodule $X$
over a semiring $R$, we say that a partial order $\leq$ is compatible with the module operations,
if $r x^1 \leq rx^2$, $x+x^1 \leq x+x^2$ for all $x \in X, r\in R$, whenever $x^1\leq x^2$.
}
\end{definition}

%{\color{red}
\begin{proposition} \label{Pro:divtwotim}{\rm  \cite[Proposition 2.2]{LM25}}   {\it Let $U$ and $V$ be two semimodules over a semiring $R$. Then $U$ and $V$ are (uniquely) divisible, if and only if $U\times V$ is (uniquely) divisible.
}
\end{proposition}

Let $(X,+)$ be a semigroup and $C, D$ be subsets of $X$.
The general Minkowski sum  of sets in $X$ is given by $C+D:=\{c+d: c\in C, d\in D\}$.
If $D:=\{d\}$, then we set $C+D:=C+d$. If $(X,+)$ is a group, then we denote by $c-d:=c+(-d)$ for $c,d\in X$ and $C-D:=C+(-D):=\{c-d: c\in C, d\in D\}$.

Given  a semimodule  $Y$ over the semiring $R$ where $Y$ is equipped with a partial order $\leq$,
we may consider the maximal and minimal elements in $Y$, $\pm\infty$, where we define
 on $Y\cup \{\pm\infty\}$  the  operations given on $Y$ to which we add the following
one: $x+(\pm\infty)=(\pm\infty)+x:=\pm\infty, \forall x\in Y$; $r\cdot (\pm\infty):=\pm\infty, \forall r\in R\backslash \{0\}$; $-\infty \leq x\leq +\infty, \forall x\in Y$; $+\infty+(-\infty):=+\infty$;
$0\cdot (\pm\infty)=0$.

Throughout the paper, we always refer to additive groups and semigroups, unless differently  specified. Moreover, without other specification we always denote by  $\pm\infty$ the maximal and minimal
elements  of corresponding semimodules over semirings, respectively.

\section{Convex sets and cones in semigroups and groups}
\setcounter{equation}{0}
In this section, we analyse in details convex sets and cones in semimodules over semirings.
\begin{definition}\label{Def:3.1} {\rm (Convex set and convex hull in a semimodule)} {\rm \cite[Definition 3.1]{LM24}} {\it  Let $X$ be a semimodule over a
semiring $R$. A set $A\subseteq X$ is called a convex set, if for any choice of $x\in X, n\in \mathbb{N}, r, r^i \in R\backslash\{0\}$ and $x^i\in A (i=1,\cdots,n)$
such that
$
rx=\sum_{i=1}^nr^i x^i,\,\, r=\sum_{i=1}^nr^i,
$
one has $x\in A$. The convex hull of $A$, say ${\rm conv} A$, is given by ${\rm conv} A:=\cap \{B: A\subseteq B \,{and}\, B\,\,{is\ convex}\}$.
}
\end{definition}

If $X$ is a semimodule over a semiring $R$ and $C\subseteq D\subseteq X$, then it is easy to see that ${\rm conv} C\subseteq {\rm conv} D$.

In the following, we introduce cones and convex cones in semimodules over semirings.

\begin{definition}\label{Def:3.3} {\rm (Cone in a semimodule)} {\rm  \cite[Definition 3.12]{LM24}}  {\it  Let $X$ be a semimodule over a semiring $R$. A set $A \subseteq X$ is called a cone,
 if for any choice $x' \in X$, $r,r' \in R \setminus \{0\}$, $x\in A$ such that $r'x'=rx$, one has $x'\in A$.
The set  $R\backslash\{0\} A:=\{x'\in X: r'x'=r x, x\in A, r, r'\in R\backslash\{0\}\}$  is called  the cone generated by $A$.}
\end{definition}
 From Definition \ref{Def:3.3} it is easy to see that $A$ is a cone if and only if $R\backslash\{0\} A=A$.

\begin{definition}\label{Def:3.4}  {\rm  \cite[Definition 3.13]{LM24}} {\rm (Convex cone in a semimodule)} %\cite[Definition 8]{BG18}}
 {\it Let $X$ be a semimodule over a
semiring $R$. A set $A\subseteq X$ is said to be a convex cone, if for any choice $x\in X$, $n\in \mathbb{N}$,  $r, r^i\in R\backslash\{0\}$ and $x^i\in A (i=1,\cdots,n)$ such that
\begin{equation}\label{cone1}
rx=\sum_{i=1}^nr^i x^i,\quad  \sum_{i=1}^nr^i\ne 0,
\end{equation}
one has $x \in A$.
The set of the points $x$ such that \eqref{cone1} holds is called  the convex cone generated by $A$ and is denoted by  ${\rm cone}\, A$.
}
\end{definition}

Let $X$ be a semimodule over a semiring $R$ and $A\subseteq X$. Then it is clear that $A\subseteq R\backslash\{0\} A\subseteq {\rm cone}\, A$.

\begin{proposition} \label{Pro:3.1} {\rm \cite[Proposition 3.4]{LM24}}  {\it  Let $X$ be a semimodule over a
semiring $R$. Then $A$ is convex if and only if $A={\rm conv} A$.
}
\end{proposition}

In case $X$ is a monoid the following characterization of the convex hull has been proved in \cite{BG18}.

\begin{proposition} \label{Pro:3.2}   {\rm \cite[ Proposition 1]{BG18}}  {\it  Let $X$ be a monoid and $A\subseteq X$. Then
\begin{eqnarray*}
{\rm conv} A=\left\{ x\in X: mx=\sum_{i=1}^n m^i x^i, \ x^i\in A, \ n, m^i\in\mathbb{N}, \ m=\sum_{i=1}^n m^i\right\}.
\end{eqnarray*}
}
\end{proposition}

Clearly, for a monoid $X$ we always have ${\rm conv} X=X$, that is, $X$ is convex.

\begin{proposition} \label{Pro:3.3}   {\it  Let $C$ be a convex subset of a commutative group $X$ with $0\in C\not=X$ and let $x^0\not\in C$. Then $x^0$ has infinite order.
}
\end{proposition}
{\bf Proof.} Suppose to the contrary, that $x^0$ has finite order, i.e., there exists $ m\in\mathbb{N}$ such that $ m x^0=0$.
Since $0\in C$ and $C$ is convex, there are $r, n, r^i\in \mathbb{N} (i=1,\cdots,n)$, with $r= \sum_{i=1}^nr^i$, and $x^i\in C (i=1,\cdots,n)$ such that
$r\cdot 0=\sum_{i=1}^nr^i x^i$.
It follows that $mr= \sum_{i=1}^n mr^i$ and
$
mrx^0=r0= 0= mr\cdot 0=\sum_{i=1}^n mr^i x^i,
$
which implies  that $x^0\in C$, a contradiction. This proves  that $x^0$ has infinite order. This completes the proof. \cvd

%{\color{red}
%\begin{proposition} \label{Pro:singletonconv} {\rm  \cite[Proposition 3.3]{LM25}}   {\it  Let $X$ be a commutative group and $A:=\{x^0\}$, where $x^0\in X$. If for each $x\in X\backslash \{0\}$, the cyclic subgroup $\langle x %\rangle$ has infinite order (especially, $X$ is uniquely divisible), then $A$ is convex and ${\rm conv}A:=\{x^0\}$.
%}
%\end{proposition}
%}

\begin{proposition} \label{Pro:3.4}  {\it  Let $A$ be a  subset of a commutative group  $X$  and $x^0\in X$. Then $A$ is convex if and only if $x^0+A$ (or equivalently, $x^0-A$) is convex.
}
\end{proposition}
{\bf Proof.} From {\rm \cite [Lemma 3.4]{LM23}} we have $A$ is convex if and only if $x^0+A$ is convex. From the definition of convexity, it is easy to check that $A$ is convex if and only if $-A$ is convex. Now, the previous conclusion yields that $x^0+A$ is convex if and only if $x^0-A$ is convex.
\cvd

\begin{proposition} \label{Pro:3.8} {\rm  \cite[Propositions 3.3 and  3.18 (iii)]{LM24}}  {\it  Let $X$ be a semimodule over a semiring $R$ and $A, B$ be convex subsets (resp., cones) of $X$. If $X$ is uniquely divisible, then $A+B$ is a convex subset (resp., cone) of $X$.
}
\end{proposition}

\begin{proposition} \label{Pro:3.4+} {\it  Let $A, B$ be convex subsets of a group $X$. If $X$ is uniquely divisible, then  $A\pm B$ are convex.
}
\end{proposition}
{\bf Proof.} Since $B$ is convex,  from Proposition \ref{Pro:3.4} we have $-B=0-B$ is convex. If $X$ is uniquely divisible, then the conclusion follows immediately from Proposition \ref{Pro:3.8}.
\cvd

Given a subset $A$ of a monoid $X$ and $m\in \mathbb{N}$. Let
\begin{eqnarray}\label{3.1}
m A&:=& \{ x\in X: x=\sum_{i=1}^m x^i, x^i\in A\}=\{ x\in X: x=\sum_{i=1}^n m^ix^i, x^i\in A, \sum_{i=1}^n m^i=m,  m^i, n\in \mathbb{N}\}.
\end{eqnarray}

We denote by $0A:=\{0\}$. Note that \eqref{3.1} can be generalized in the following in the case of convex sets.

\begin{definition}\label{Def:3.2} {\rm (Rational dilation of set)}   {\it Let $C$ be a subset of a monoid $X$ and let $q\in \mathbb{Q}_+\backslash \{0\}$
be a reduced fraction. Define
\begin{eqnarray}\label{3.2}
q C:= \{ x\in X: l x=\sum_{i=1}^n m^ix^i, x^i\in C, \sum_{i=1}^n m^i=m, \frac{m}{l}=q, l, n, m^i \in \mathbb{N}\}.
\end{eqnarray}
}
\end{definition}
In {\rm \cite[Definition 2.6]{BG16}}, $C$ is assumed to be a convex subset of a monoid $X$.
It is easy to see that $1 C\supseteq C$.
Moreover, if $C$ is not convex, then the equality $1 C= C$ does not necessarily hold,
in fact we note that $1 C={\rm conv}\, C$.

\begin{proposition} \label{Pro:3.5}   {\rm  \cite[Proposition 2.3]{BG16}}  {\it
Let $X$ be a monoid, $C\subseteq X$ and let  $q^1, q^2\in \mathbb{Q}_+$ be reduced fractions. Assume that  $0\in C$ and $C$ is convex. If $q^1\leq q^2$, then $q^1C\subseteq q^2C$.
}
\end{proposition}

\begin{proposition} \label{Pro:3.5+}  {\rm  \cite[Proposition 3.6]{LM25}}  {\it Let $C$ be a subset of a monoid $X$ and $q\in \mathbb{Q}_+$ be a reduced fraction.  If $0\in C$, then $0\in q C$.
}
\end{proposition}

\begin{proposition} \label{Pro:3.6}    {\rm  \cite[Proposition 3.15]{LM24}} {\it Let $X$ be a semimodule over a
semiring $R$ and $C\subseteq X$. Then the following statements hold:
\begin{itemize}
\item[\rm(i)]  $C$ is a convex cone if and only if $C={\rm cone}\, C$, which implies that $C$ is both a cone and convex.

\item[\rm(ii)] Let $X$ be  a monoid. Then ${\rm cone}\, C=\displaystyle\cup_{q\in \mathbb Q_+\backslash\{0\} }qC$, which implies that $C=\displaystyle\cup_{q\in \mathbb Q_+\backslash\{0\} }qC$ if $C$ is a convex cone.

\item[\rm(iii)] Let $X$ be a monoid. If $C$ is a convex cone, then $C$ is a cone and $C+C\subseteq C$, and so $C+C=C$ if $0\in C$.
Conversely, if $C$ is a cone and $C+C=C$, then $C$ is a convex cone.

\item[\rm(iv)] Then $C$ and $D$ are convex cones if and only if $C\times D$ is a convex cone.
\end{itemize}
}
\end{proposition}

\begin{proposition} \label{Pro:3.7}  {\rm  \cite[Proposition 3.10]{LM25}}  {\it Let $X$ be a divisible monoid and $C\subseteq X$. If $C$ is a convex cone, then $C=qC$ for any reduced fraction $q\in \mathbb Q_+\backslash\{0\}$.
}
\end{proposition}

\begin{proposition} \label{Pro:3.10} {\it  Let $A$ be a  (resp., convex) cone of a monoid  $X$. Then  $a\in A$ if and only if $ra\in A$, where $r\in \mathbb{N}$.
}
\end{proposition}
{\bf Proof.} Necessity. It follows immediately from  Definition \ref{Def:3.3} (resp., Proposition \ref{Pro:3.6} (iii)).

Sufficiency. Let $ra\in A$ and $r\in \mathbb{N}$. Then there is $b\in A$ such that $ra=b$, which yields that $a\in A$ in view of  Definition \ref{Def:3.3} (resp., Proposition \ref{Pro:3.6} (i)).
\cvd

In the following proposition, we remove the assumption that the divisible commutative group $X$ is unique in {\rm  \cite[Proposition 3.18 (iii)]{LM24}}.

\begin{proposition} \label{Pro:3.9} {\it  Let $A$ and $B$ be (resp., convex) cones of a commutative group $X$. If $X$ is divisible, then $A\pm B$ are (resp., convex) cones.
}
\end{proposition}
{\bf Proof.} It is easy to check that for a subset $C$ of $X$, $C$ is a (resp., convex) cone if and only if $-C$ is a (resp., convex) cone. It suffices to show $A+B$ is a (resp., convex) cone. Similarly, we can show $A-B$  is a (resp., convex) cone.

Let  $z'\in X$, $r', r\in \mathbb{N}$ and $z\in A+B$ such that
$r'z'=r z$. Then, there are $x, y\in X$ such that $z=x+y$ and so $r'z'=rx+ry$. Since $X$ is divisible, there are $x', y'\in X$ such that $r'x'=rx$ and $r'y'=ry$ and from Proposition \ref{Pro:3.10} it follows that $r'(z'-y')=r'x'\in A$. Again  from Proposition \ref{Pro:3.10}, one has $z'-y' \in A$, i.e., $z' \in A+y'\subseteq A+B$, which yields that $A+B$ is a cone.

Let  $z\in X$, $n, r, r^i\in \mathbb{N}$, $z^i\in A+B (i=1,\cdots,n)$ such that $rz=\sum_{i=1}^nr^i z^i$.
Then, there are $x^i\in A$ and $y^i\in A (i=1,\cdots,n)$  such that $z^i=x^i+y^i (i=1,\cdots, n)$ and so
$rz=\sum_{i=1}^nr^i (x^i+y^i)$. Since $X$ is divisible, there are $x, y\in X$ such that $rx= \sum_{i=1}^nr^ix^i$ and $ry= \sum_{i=1}^nr^iy^i$, which yields that $x\in A$ and $y\in B$, since $A$ and $B$ are convex cones.
As a consequence, $rz=rx+ry$ and from Proposition \ref{Pro:3.6} (iii) it follows that $r(z-y)=rx\in A$. This, along with Proposition \ref{Pro:3.10}, yields that $z-y \in A$ and so $z \in A+y\subseteq A+B$.
This proves that $A+B$  is a convex cone.
\cvd

\begin{proposition} \label{Pro:conesubgroup} {\rm  \cite[Lemma 4.5]{LM25}} {\it Let $A$ be a subgroup of a  group $X$. Then ${\rm cone}\,A$ is a convex subgroup of $X$.
}
\end{proposition}

\begin{proposition} \label{Pro:CPO} \rm  \cite[Proposition 3.1]{LM23}    {\it  Assume that $Y$ is a commutative group, $D\subseteq Y$ is a pointed (i.e., $(-D)\cap D=\{0\}$) convex cone. Let $x,y\in Y$. Define $y\leq x \Longleftrightarrow x-y\in D$. Then $\leq$ is a compatible partial order.
}
\end{proposition}

\begin{definition} \label{Def:3.5}   {\it  Assume that $C$ is a subset of  a group $X$.
\begin{itemize}
\item[\rm(i)] {\rm (Core)} {\rm \cite[Definition 26]{BG18}}  The core of $C$ is defined by
$
{\rm core}\, C:=\{x\in X: \forall h\in X, \exists n\in \mathbb{N}, \exists g\in X, ng=h, x+g\in C\}.
$

\item[\rm(ii)]  {\rm (Quasi core)} \rm \cite[Definition 4.21]{LM24}  The quasi core of  $C$ is defined by
$
{\rm qcore}\,C:=\{x\in C:{\rm cone}\,(C-x)=X\}.
$

\item[\rm(iii)] {\rm (Pseudo relative core)} \rm \cite[Definition 4.1]{LM25} The pseudo relative core of  $C$ is defined by
\begin{eqnarray*}
{\rm prcore}\,C&:=&\{x\in C:{\rm cone}\,(C-x)\, {\rm is\,\, a\,\, subgroup\,\, of}\, X\}\\
&:=&\{x\in C:{\rm cone}\,(C-x)\, {\rm is\,\, a\,\, convex\,\,subgroup\,\, of}\, X\}.
 \end{eqnarray*}
%}
\end{itemize}
}
\end{definition}
The second equality in  Definition \ref{Def:3.5} (iii) follows immediately from Proposition \ref{Pro:3.6} (i).
From \cite[Lemma 4.22 (i)]{LM24}, we have ${\rm core}\, C\subseteq {\rm qcore}\,C\subseteq C$ if $X$ is a commutative  group. Clearly, ${\rm qcore}\,C\subseteq {\rm prcore}\,C\subseteq C$ and ${\rm qcore}\, C\subseteq {\rm qcore}\, D$,  if $C$ and $D$ are subsets of $X$ such that $C\subseteq D$. Several interesting properties and examples of quasi (relative) core of a subset of a group were investigated by Li and Mastroeni in \cite{LM24, LM25}.

\begin{proposition} \label{Pro:eqqqrcore} \cite[Lemma 4.25]{LM24}  \cite[Proposition 4.9]{LM25}  {\it Assume that $X$ is a commutative group, $\bar x\in X$ and  $\emptyset \neq C\subseteq  X$. Then $ {\rm core}\,(\bar x +C)= \bar x + {\rm core}\,C$, $ {\rm qcore}\,(\bar x +C)= \bar x + {\rm qcore}\,C$ and  $ {\rm prcore}\,(\bar x +C)= \bar x + {\rm prcore}\,C$.
}
\end{proposition}

%{\color{red}
%{\bf Proof.} Observe that $B$ is a  (resp., convex) cone if and only if $-B$ is a (resp., convex) cone. %Assume that $X$ is divisible.
%It suffices to show $A+B$ ia a convex cone. Similarly, we can show $A-B$ ia a convex cone. Let  $z\in X$, $n, r, r^1, \cdots, r^n\in \mathbb{N}$, $x^1, \cdots, x^n\in A$ and $y^1, \cdots, y^n\in A$  such that
%$rz=\sum_{i=1}^nr^i (x^i+y^i)$.

%By Proposition \ref{Pro:3.10}, we have that $r^ix^i\in A$ and $r^iy^i\in B$. From Proposition \ref{Pro:3.6} (iii) it follows that
%$$x= \sum_{i=1}^nr^ix^i\in A\quad and \quad y= \sum_{i=1}^n r^i y^i\in B$$

 %$$rz=\sum_{i=1}^nr^i (x^i+y^i)=x+y\in A+B.$$

%As a consequence, $rz=rx+ry$ and from Proposition \ref{Pro:3.6} (iii) it follows that $r(z-y)=rx\in A$. This, along with Proposition \ref{Pro:3.10}, yields that $z-y \in A$ and so $z \in A+y\subseteq A+B$.
%This proves that $A+B$ ia a convex cone.
%\cvd
%}

\section{Convex functions in semigroups and groups}
\setcounter{equation}{0}
In this section, we study convex functions in groups and semigroups.

\begin{definition}\label{Def:4.1} {\rm  (Epigraph, effective domain and lower level set) }  {\it Suppose that $X, Y$ are semimodules over a semiring $R$ and that $Y$ is endowed with a partial order $\leq$. Let $f:X\rightarrow Y\cup \{\pm\infty\}$.
\begin{itemize}
\item[\rm(i)] The epigraph of $f$, say ${\rm epi}\,f$, is
given by
${\rm epi}\,f:=\{(x,y)\in X\times Y: f(x)\leq y\}$.

\item[\rm(ii)] The effective domain of $f$, say ${\rm dom}\,f$, is
given by
${\rm dom}\,f:=\{x\in X: f(x)<+ \infty\}=\{x\in X: \exists y\in X, (x,y)\in {\rm epi}\,f\}$.

\item[\rm(iii)] The lower level set of $f$, say ${\rm lev}_{\leq y}\,f$, is
given by ${\rm lev}_{\leq y}\,f:=\{x\in X: f(x)\leq y\}$, where $y\in Y$.
\end{itemize}
}
\end{definition}

If $f:X\rightarrow Y\cup \{+\infty\}$, then previous definition collapses to {\rm \cite[Definition 4.1]{LM24}}.

\begin{definition}\label{Def:4.2} {\rm  (Convex and concave function) } {\it Suppose that $X, Y$ are  semimodules over a semiring $R$ and that $Y$ is endowed with a partial order $\leq$. A function $f: X\to Y\cup \{\pm\infty\}$
is said to be convex (resp., concave), if for any choice of $x\in X$, $n\in \mathbb{N}$, $r, r^i \in R\backslash \{0\}$ and
$x^i \in X (i=1,\cdots,n)$ such that
$
rx=\sum_{i=1}^nr^i x^i,\,\, r=\sum_{i=1}^nr^i,
$
we have
$
rf(x)\leq \sum_{i=1}^nr^i f(x^i) \quad (resp., rf(x)\geq \sum_{i=1}^nr^i f(x^i)).
$
}
\end{definition}

If $f:X\rightarrow Y\cup \{+\infty\}$, then previous definition collapses to {\rm \cite[Definition 4.2]{LM24}}.
The definition of the convex function $f: X\to Y$ was given in \cite[Definition 16]{BG18}  without the assumption $r\ne 0$.

%It is easy to see that $f: X\to Y\cup \{+\infty\}$ is convex if and only if $f: {\rm dom}f\to Y$ is convex.
%In the following, without other specification, we always assume that (A1) for every $r\in R\setminus \{0\}, y\in Y\cup \{+\infty\}$, $r y=+\infty %\Longleftrightarrow y=+\infty$.

%{\color{red} Actually, in \cite[Definition 16]{BG18} it is considered $f:X\to Y$.
%$As in standard convex analysis, maybe we could  say in the previous definition that $f: X\to Y\cup \{+\infty\}$ is convex if for every $n\in \mathbb{N}$, every $r_1, \cdots, r_n \in R\backslash \{0\}$ and every
%$x_1,\cdots, x_n \in dom\, f$ such that
%\begin{eqnarray*}
%rx=\sum_{i=1}^nr_i x_i,\,\, r=\sum_{i=1}^nr_i,
%\end{eqnarray*}
%we have
%\begin{eqnarray*}
%rf(x)\leq \sum_{i=1}^nr_i f(x_i) .
%\end{eqnarray*}
%}

The proofs of (ii)-(iv) in the following proposition are similar to that in  {\rm \cite[Proposition 4.3]{LM24}}. For completeness, we include it here.

\begin{proposition} \label{Pro:4.1}   {\it Suppose that $X, Y$ are semimodules over a semiring $R$  and that $Y$ is endowed with a compatible partial order $\leq$. Let $f:X\rightarrow Y\cup \{\pm\infty\}$. Then the following statements hold:
\begin{itemize}
\item[\rm(i)] $f:X\rightarrow Y\cup \{+\infty\}$ is convex if and only if $f:{\rm dom}\,f \rightarrow Y$ is convex.

\item[\rm(ii)]  If $f$ is convex, then ${\rm dom}\,f$ is convex.

\item[\rm(iii)] Assume that $(Y, \leq)$ satisfies that for every $r\in R\setminus \{0\}, y^1, y^2\in Y$, $ ry^1\leq r y^2\Longrightarrow y^1\leq y^2$.  If $f$ is convex, then ${\rm epi}\,f$ and ${\rm lev}_{\leq y}\,f$ are convex, where $y\in Y$.

\item[\rm(iv)] Assume that  $Y$ is divisible and $f:X\rightarrow Y\cup \{+\infty\}$.  If ${\rm epi}\,f$ is convex, then $f$ is convex.
\end{itemize}
}
\end{proposition}
{\bf Proof.}\, (i) Necessity. It is clear.

Sufficiency. Assume that $f:{\rm dom}\,f \rightarrow Y$ is convex. Let $x\in X,  n\in \mathbb{N}, x^i\in X$ and  $r, r^i\in R\setminus \{0\} (i=1,\cdots, n)$ such that
$rx=\sum_{i=1}^nr^i x^i$ and $r=\sum_{i=1}^nr^i$. If $x^i\in {\rm dom}\,f$ for each $i=1,\cdots, n$, then it follows that $rf(x)\leq \sum_{i=1}^nr^i f(x^i)$. If $x^i\not\in {\rm dom}\,f$, i.e., $f(x^i)=+\infty$, for some $i=1,\cdots, n$, then previous inequality holds trivially. This shows that  $f:X\rightarrow Y\cup \{+\infty\}$ is convex.

(ii) Let $x\in X, n\in \mathbb{N}, x^i\in {\rm dom}\,f$ and $r, r^i\in R\setminus \{0\} (i=1,\cdots, n)$ such that
$rx=\sum_{i=1}^nr^i x^i$ and $r=\sum_{i=1}^nr^i$.
Since $x^i\in {\rm dom}\,f(i=1,\cdots, n)$, one has $f(x^i)\neq +\infty$.  Since $f$ is convex, one has
$
 rf(x)\leq \sum_{i=1}^n r^i f(x^i)\neq +\infty,
$
which yields that $f(x)\neq +\infty$, i.e., $x\in {\rm dom}\,f$. If not, that is, $f(x)= +\infty$,  then $r f(x)= +\infty$,
a contradiction.

(iii) Let $(x,y)\in X\times Y, n\in \mathbb{N}, (x^i, y^i)\in {\rm epi}\,f$ and  $r, r^i\in R\setminus \{0\} (i=1,\cdots, n)$ such that
$r(x, y)=\sum_{i=1}^nr^i (x^i, y^i)$ and $r=\sum_{i=1}^nr^i$.
Since $(x^i, y^i)\in {\rm epi}\,f (i=1,\cdots, n)$, we have $f(x^i)\leq y^i$. Since $f$ is convex, it follows from the fact that $\leq$ is a compatible partial order that
$
 rf(x)\leq \sum_{i=1}^n r^i f(x^i)\leq \sum_{i=1}^n r^i y^i= ry.
$
Now the assumption yields that $f(x)\leq y$, i.e., $(x,y)\in {\rm epi}\,f$. This proves that ${\rm epi}\,f$ is convex. Similarly, we can prove ${\rm lev}_{\leq y}\,f$ is convex.

(iv)  From (i) it suffices to show $f:{\rm dom}\,f \rightarrow Y$ is convex.
To this aim, let $x\in X,  n\in \mathbb{N}, x^i\in {\rm dom}\,f$ and  $r, r^i\in R\setminus \{0\} (i=1,\cdots, n)$ such that
$rx=\sum_{i=1}^nr^i x^i$ and $r=\sum_{i=1}^nr^i$.
Observe that  $(x^i, f(x^i))\in {\rm epi}\,f (i=1,\cdots, n)$ since $f:X\rightarrow Y\cup \{+\infty\}$.
Since  $Y$ is divisible, then  there exists $y\in Y$ such that $r y=\sum_{i=1}^nr^if(x^i)$. As a consequence,
$r(x, y)=\sum_{i=1}^nr^i (x^i, f(x^i))$.
Since ${\rm epi}\,f$ is convex, one has $(x, y)\in {\rm epi}\,f$, that is, $f(x)\leq y$ and it follows from the fact that $\leq$ is a compatible partial order that
$
 rf(x)\leq  ry= \sum_{i=1}^n r^i f(x^i).
$
This proves $f:{\rm dom}\,f \rightarrow Y$ is convex.
\cvd

\begin{definition}\label{Def:4.3} {\rm  (Minkowski functional for monoids})  {\rm \cite[Definition 2.7]{BG16}}  {\it Let $X$ be a monoid and $C\subseteq X$. The Minkowski functional
$\rho_C: X\to \mathbb{R}_+\cup \{+\infty\} $ is defined by
$
\rho_C(x):=\inf\{q\in \mathbb{Q}_+: x\in qC\}.
$
If  $\{q\in \mathbb{Q}_+: x\in qC\}=\emptyset$, then define $\rho_C(x):=+\infty$.
}
\end{definition}

\begin{remark}\label{Rm:Mink} {\it It is easy to see that if $x\in C$, then $\rho_C(x)\leq 1$. If we assume $0\in C$,
then $\rho_C(x)\geq1$, for any $x\not\in C$. Indeed, ab absurdo, suppose  that  $\rho_C(x)<1$. Then there is a reduced fraction $q\in \mathbb{Q}_+$ with $q<1$ such that
$x\in qC\subseteq C$, where the inclusion follows from Proposition {\rm\ref{Pro:3.5}}. This is a contradiction.

Form Proposition {\rm\ref{Pro:3.6} (ii)} we have $\{q\in \mathbb{Q}_+: x\in qC\}=\emptyset$ if and only if $x\not\in \displaystyle\cup_{q\in \mathbb Q_+\backslash\{0\} }qC\cup \{0\}={\rm cone}\, C\cup \{0\}$, which implies that if $x\not\in  {\rm cone}\, C\cup \{0\}$, then $\rho_C(x)=+\infty$.
}
\end{remark}

\begin{definition}\label{Def:4.4} {\rm  ($\mathbb{N}$-sublinear function) }{\rm \cite[Definitions 19 and 20]{BG18}} {\it Assume that $X, Y$ are monoids, and assume that $Y$ is equipped with a partial order $\leq$. A function $f:X\rightarrow Y\cup \{\pm\infty\}$ is said to be $\mathbb{N}$-sublinear if it is subadditive, i.e.,
$
f(x+y)\leq f(x)+f(y)
$
for any $x, y\in X$ and in addition it is positively homogeneous, i.e.,
$
f(k x)=kf(x)
$
for any $k\in \mathbb{N}\cup\{0\}$ and $x\in X$.
}
\end{definition}

\begin{remark}\label{Rm:NSunFunc} {\it In Definition \ref{Def:4.4}, we assume that $X, Y$ are monoids, but not semimodules over a semiring $R$ as in \cite[Definitions 19 and 20]{BG18}, since
$k\in \mathbb{N}\cup\{0\}$. Moreover, Borwein  and Giladi defined  $0\cdot (+\infty)=+\infty$ in \cite[Definition 20]{BG18}. In fact, if $f\equiv +\infty$, then $+\infty=f(0)=f(0 y_0)=0f(y_0)=0\cdot (+\infty)$ for any $y_0\in X$. Assume that $f\not\equiv \pm\infty$, $f(y_0)=\pm\infty$ and $-\infty<f(x_0)<+\infty$  for some $x_0, y_0\in X$. Since the $\mathbb{N}$-sublinear function $f:X\rightarrow Y\cup \{\pm\infty\}$ is  positively homogeneous, i.e., $f(k x)=kf(x)$ for any $k\in \mathbb{N}\cup\{0\}$ and $x\in X$, we have $0=0f(x_0)=f(0x_0)=f(0 y_0)=0f(y_0)=0\cdot (\pm\infty)$, as defined at the end of section 2.

The two cases, $f\equiv +\infty$ and $f\equiv -\infty$, are trivial.
In the sequel of the paper, for $f:X\rightarrow Y\cup \{\pm\infty\}$,  without other specification,  we always assume that $-\infty<f(x_0)<+\infty$  for some $x_0\in X$, which implies that $f\not\equiv \pm\infty$.
}
\end{remark}

The following proposition is an improvement of  {\rm \cite[Proposition 4.8]{LM24}} by extending $f:X\rightarrow \mathbb{R}\cup \{+\infty\}$ to $f:X\rightarrow \mathbb{R}\cup \{\pm\infty\}$ and removing the assumption that $X$ is semidivisible in (ii).

\begin{proposition} \label{Pro:4.1+0}   {\it Assume that $X$ is a monoid and $f:X\rightarrow \mathbb{R}\cup \{\pm\infty\}$. Then the following statements hold:
\begin{itemize}
\item[\rm(i)]  If $f$ is  $\mathbb{N}$-sublinear, then it is convex and positively homogeneous.

\item[\rm(ii)] If $f$ is convex and positively homogeneous, then $f$ is $\mathbb{N}$-sublinear.
\end{itemize}
}
\end{proposition}
{\bf Proof.}\, (i) Let $x\in X, n\in \mathbb{N}, x^i\in X$ and $r^i\in \mathbb{Z}_+\backslash \{0\}=\mathbb{N} (i=1,\cdots, n)$ such that
$
rx=\sum_{i=1}^nr^i x^i,\,\, r=\sum_{i=1}^nr^i.
$
Since $f$ is  $\mathbb{N}$-sublinear, it is positively homogeneous and subadditive and as a consequence,
\begin{eqnarray*}
rf(x)=f(rx)=f(\sum_{i=1}^nr^i x^i)\leq \sum_{i=1}^nf(r^i x^i)=\sum_{i=1}^n r^i f(x^i),
\end{eqnarray*}
which implies that $f$ is convex.

(ii)  Let $x, y\in X$ and set $z:=x+y$. Then $2z=2(x+y)=(2x)+(2y)$. Since $f$ is convex and positively homogeneous, it follows that
$
2f(z)\leq f(2x)+f(2y)=2f(x)+2f(y)
$
and so
$f(z)\leq f(x)+f(y)$,
yielding that $f$ is subadditive and so $\mathbb{N}$-sublinear.
\cvd

 The second statement of following lemma is a slight modification of \cite[Lemma 4.22]{LM24}.

\begin{proposition} \label{Pro:4.2} {\it Let  $X$ be  a monoid and $C\subseteq X$. Then the following statements hold:
 \begin{itemize}
\item[\rm(i)] $\rho_C:X\to \mathbb{R}_+\cup \{+\infty\}$ is $\mathbb{N}$-sublinear and so convex and positively homogeneous.

\item[\rm(ii)]  $0\leq\rho_C(x) < +\infty$ for any $x\in {\rm cone}\,C$.

\item[\rm(iii)] If $C$ is a convex cone and $X$ is divisible, then $\rho_C(x)=0$ for any $x\in {\rm cone}\,C$.

\end{itemize}
}
\end{proposition}
{\bf Proof.}  (i) It follows from  \cite[Proposition 4.9 (i)]{LM24} and  Proposition \ref{Pro:4.1+0} (i) (see also, \cite[Proposition 2.4]{BG16}, in fact, in which the convexity of $C$ can be removed).

(ii) By Proposition  \ref{Pro:3.6} (ii), we have ${\rm cone}\, C = \cup_{q\in\mathbb{Q}_+\setminus\{0\}} q\,C$ and thus for every $x\in {\rm cone}\, C$, there exists $q^x\in \mathbb{Q}_+\setminus\{0\}$ such that $x\in q^x\,C$. This implies  $0\leq \rho_C(x) \leq q^x< +\infty$ for any $x\in {\rm cone}\,C$.

(iii) If $C$ is a convex cone and $X$ is divisible, then from Propositions \ref{Pro:3.7} and \ref{Pro:3.6} (i)  we have $C=qC={\rm cone}\,C$ for any reduced fraction $q\in \mathbb Q_+\backslash\{0\}$.
This allows that $\rho_C(x)=0$ for any $x\in {\rm cone}\,C$.
\cvd

 \begin{definition}\label{Def:4.5} {\rm (Generalised affine function) }{\rm  \cite[Definition 18]{BG18}} {\it Assume that $X$ and $Y$ are semimodules over a semiring $R$  and assume that $Y$ is equipped with a partial order $\leq$. The function $f:X\to Y \cup \{\pm\infty\}$ (resp., $f:X\to Y$) is said to be generalised   affine (resp., affine) if it is both convex and concave.}
\end{definition}

 According to Definition \ref{Def:4.2} we have the following characterization of a generalised affine function.

\begin{proposition}\label{Pro:4.2a}  {\it Suppose that $X, Y$ are  semimodules over a semiring $R$, and assume that $Y$ is equipped with a partial order $\leq$ and that  $f: X\to Y\cup \{\pm\infty\}$
is  generalised affine. Then,  $f: X\to Y\cup \{\pm\infty\}$
is  generalised affine if and only if for any choice of $n\in \mathbb{N}$, $r, r^i \in R\backslash \{0\}$ and
$x, x^i \in X (i=1,\cdots,n)$ such that
$
rx=\sum_{i=1}^nr^i x^i,\,\, r=\sum_{i=1}^nr^i,
$
we have
$
rf(x) = \sum_{i=1}^nr^i f(x^i) .
$
}
\end{proposition}

In \cite[Proposition 4]{BG18}, the group $X$ should be assumed to be commutative from its proof.

\begin{proposition} \label{Pro:4.3}{\rm  \cite[Proposition 4]{BG18}} {\it Let $X$ and $Y$ be commutative groups, and $Y$ be equipped with a partial order $\leq$ and $f:X\to Y\cup \{\pm\infty\}$  be generalised affine. Then
either  $f$  is everywhere finite, or $f\equiv+\infty$, or $f\equiv -\infty$, or $f$ attains both values $+\infty$ and $-\infty$.
}
\end{proposition}

The following is a correction of Proposition \ref{Pro:4.3}.

\begin{proposition} \label{Pro:refinementof4.3}  {\it Let $X$ and $Y$ be commutative groups, and $Y$ be equipped with a partial order $\leq$ and $f:X\rightarrow Y\cup \{\pm\infty\}$ be a generalised affine function. Then either $f:X\to Y$, or $f\equiv+\infty$, or $f\equiv-\infty$.
or $f$ can not attain both values $+\infty$ and $-\infty$.
}
\end{proposition}
{\bf Proof.}  It suffices to show it is impossible that $f$ attains both values $+\infty$ and $-\infty$.  Suppose to the contrary that there are $x^0, y^0\in X$ such that $f(x^0)=+\infty$ and $f(y^0)=-\infty$. Since $2\cdot 0=x^0+(-x^0)$ and  $f$ is generalised affine, from Proposition \ref{Pro:4.2a} we have $2f(0)=f(x^0)+f(-x^0)=+\infty$   since $+\infty+t=+\infty$ for any $t\in Y\cup \{\pm\infty\}$.  This yields that $f(0)=+\infty$.
Since $2y^0=(2y^0)+0$ and  $f$ is generalised affine, again from Proposition \ref{Pro:4.2a} we have $2f(y^0)=f(2y^0)+f(0)=+\infty$ and so $f(y^0)=+\infty$, a contradiction with $f(y^0)=-\infty$.
This completes the proof.
\cvd

\begin{proposition} \label{Pro:4.4}  {\it Assume that $X$ is a monoid and $f:X\rightarrow \mathbb{R}\cup \{\pm\infty\}$. Then the following statements hold:
\begin{itemize}
\item[\rm(i)]  If $\pm f$ are  $\mathbb{N}$-sublinear, then $f$ is a generalised affine function and positively homogeneous.

\item[\rm(ii)] Assume that $X$ is a commutative group. If  $f$ is generalised affine with $f(0)=0$, then $f:X\rightarrow \mathbb{R}$ and $f$ is positively homogeneous.
%\item[\rm(iii)] If $X$ is semidivisible, and $f:X\rightarrow \mathbb{R}$ is an  affine function, then $\pm f$ are $\mathbb{N}$-sublinear.
\end{itemize}
}
\end{proposition}
{\bf Proof.} (i) It follow immediately from Propositions \ref{Pro:4.2a} and \ref{Pro:4.1+0} (i).

(ii) Assume that $f$ is generalised affine with $f(0)=0$. For every $y\in X$, we have $2\cdot 0=y+(-y)$ and so, taking into account Proposition \ref{Pro:4.2a}, it follows
that
$0=2f(0)=f(y)+f(-y)$,
which implies that $f(y)\in \mathbb{R}$ for any $y\in \mathbb{R}$,  i.e., $f:X\rightarrow \mathbb{R}$, taking into account that $+\infty+(-\infty)=+\infty$.

Let $x\in X$ and $k\in \mathbb{N}\cup\{0\}$. If $k\in \mathbb{N}$, then from the fact that $f:X\rightarrow \mathbb{R}$ is affine with $f(0)=0$ we have  $kx=(kx)+\overbrace{0+\cdots+0}^{k-1}$ and so  $kf(x)=f(kx)+\overbrace{f(0)+\cdots+f(0)}^{k-1}=f(kx)$.
If $k=0$, then from  $f:X\rightarrow \mathbb{R}$ with $f(0)=0$ it follows  that
$f(0x)=f(0)=0=0f(x)$. Both cases imply that $f$ is positively homogeneous. This completes the proof.
%(iii) follows immediately from (ii) and Proposition \ref{Pro:4.1+0} (i).
\cvd

\begin{definition}\label{Def:additive} {\rm  (Additive function) } {\it Assume that $X, Y$ are semimodules over a semiring $R$ and assume that $Y$ is equipped with a partial order $\leq$. A function $a:X\rightarrow Y\cup \{\pm\infty\}$ is said to be additive  if $a(x+y)=a(x)+a(y)$ for all $x,y\in Y$.
}
\end{definition}

Given two monoids $X$ and $Y$, let $L(X, Y)$ be the collection of all additive maps  $a:X\to Y$.
If $Y:=Z$ and  $X, Z$ are groups, then  $L(X, Z):=\{a: X\to Z: a\,\, {{\rm is}\,\, {\rm additive}}\}$ is the additive
dual group of a group $X$, denoted by $X^*$, which is an additive group with the addition being point-wise addition of elements
in $Z$.  Clearly, if $Z$ is a commutative group, then so is $Z^*$.
In particular, if $Z:=\mathbb{R}$, then we denote by $D^*:=\{a\in X^*: a(x)\geq 0, \forall x\in D\}$ the dual (positive polar) cone of a nonempty subset $D$ of $X$. It is easy to check that $D^*$ is a convex cone of $X^*$. The dual (positive polar) cone of $D^*$, say $D^{**}$, is defined by $D^{**}:=\{x\in X: a(x)\geq 0, \forall a\in D^*\}$. Let $\lambda:=(\alpha, \beta)\in (X\times Y)^*:=X^*\times Y^*$, where $X$ and $Y$ are two groups. Define $\lambda(u,v):=(\alpha, \beta)(u,v):=\alpha(u)+\beta (v)$ for $(u,v)\in X\times Y$.

\begin{proposition} \label{Pro:4.4+} {\rm \cite[Lemma 3.6]{LM23}} {\it Let $A$ be a nonempty subset of a monoid $X$ and $f:X\to \mathbb{R}\cup\{\pm \infty\}$  be an additive function. Then
$\sup_{x\in A}f(x) = \sup_{x\in {\rm conv}A}f(x) \,\, {\rm and}\,\, \inf_{x\in A}f(x) = \inf_{x\in {\rm conv}A}f(x)$.
}
\end{proposition}

\begin{proposition} \label{Pro:additiveope} {\rm\cite[Proposition 3.7 (iii)]{LM23}} {\it Let $A$ be a nonempty subset of a commutative group $X$.  Then
\begin{eqnarray*}
A^*=({\rm conv}\,A)^*=({\rm cone}\,A)^*=({\rm conv}\,({\rm cone}\,A))^*=({\rm cone}\,({\rm conv}\,A))^*.
\end{eqnarray*}
}
\end{proposition}

\begin{proposition} \label{Pro:additive char}  {\it Let $X$ be a commutative group and $f:X\rightarrow \mathbb{R} \cup \{\pm\infty\}$ be an additive function. Then the following statements are true:
\begin{itemize}
\item[\rm(i)]  $f$ is generalised affine such that either $f:X\to \mathbb{R}$ with $f(0)=0$, or $f\equiv+\infty$, or $f\equiv -\infty$.

\item[\rm(ii)] $f:X\to \mathbb{R}$ if and only if $f(0)=0$, which yields that $f(-x)=-f(x)$ for any $x\in X$.

\item[\rm(iii)] $f\equiv+\infty$ if and only if $f(0)=+\infty$.

\item[\rm(iv)] $f\equiv-\infty$ if and only if $f:X\rightarrow \mathbb{R} \cup \{-\infty\}$ with $f(0)=-\infty$.

\item[\rm(v)] ${\rm epi}\,f$ and ${\rm dom}\,f$ are convex cones.
\end{itemize}
}
\end{proposition}
{\bf Proof.} (i) Let $n, r, r^i \in \mathbb{N}$ and
$x, x^i \in X (i=1,\cdots,n)$ such that
$rx=\sum_{i=1}^nr^i x^i,\,\, r=\sum_{i=1}^nr^i$.
Since $f$ is  additive, it follows that
\begin{eqnarray*}
rf(x)=f(rx)=f(\sum_{i=1}^nr^ix^i)=\sum_{i=1}^nf(r^ix^i)=\sum_{i=1}^n r^if(x^i),
\end{eqnarray*}
which yields that $f$ is generalised affine. Now, from Proposition \ref{Pro:refinementof4.3} we have either $f:X\to \mathbb{R}$, or $f\equiv+\infty$, or $f\equiv -\infty$, or $f$ can not attain both values $+\infty$ and $-\infty$. If $f(y_0)=\pm\infty$ and $-\infty<f(x_0)<+\infty$  for some $x_0\in X$, then $f(x_0)=f(x_0+0)=f(x_0)+f(0)$ and so $f(0)=0$, since $f$ is additive. Also, we have $0=f(0)=f(x-x)=f(x)+f(-x)$ and so $f(x), f(-x)\in \mathbb{R}$ for any $x\in X$, a contradiction with $f(y_0)=\pm\infty$. This proves that $f$ is generalised affine such that either $f:X\to \mathbb{R}$ with $f(0)=0$, or $f\equiv+\infty$, or $f\equiv -\infty$.

(ii) Necessity. Since $f:X\to \mathbb{R}$ is additive, $f(0)=f(0+0)=f(0)+f(0)$ and so $f(0)=0$.

Sufficiency. Since $f(0)=0$, one has $0=f(0)=f(x-x)=f(x)+f(-x)$ and so $f(x), f(-x)\in \mathbb{R}$ for any $x\in X$. This yields that $f:X\to \mathbb{R}$ and  $f(-x)=-f(x)$ for any $x\in X$.

(iii) Necessity. It is clear.

Sufficiency. Assume that $f(0)=+\infty$. Since $f$ is  additive,  $f(x)=f(x+0)=f(x)+f(0)=+\infty$ for any $x\in X$, i.e., $f\equiv+\infty$.

(iv) Necessity. It is clear.

Sufficiency. Assume that $f:X\rightarrow \mathbb{R} \cup \{-\infty\}$ with $f(0)=-\infty$. Observe that $(-\infty)+t=-\infty$ for any $t\in \mathbb{R} \cup \{-\infty\}$.
Since $f$ is additive,  $f(x)=f(x+0)=f(x)+f(0)=-\infty$ for any $x\in X$, i.e., $f\equiv-\infty$.

(v) It suffices to show ${\rm epi}\,f$ is a convex cone. Similarly, we can show  ${\rm dom}\,f$ is a convex cone.
If ${\rm epi}\,f=\emptyset$, then it is trivial. Assume that ${\rm epi}\,f\neq\emptyset$.
Let $n, r, r^i \in \mathbb{N}$ and
$(x,\alpha)\in X\times \mathbb{R}, (x^i, \alpha^i) \in {\rm epi}\,f (i=1\cdots, n)$ such that
$r(x,\alpha)=\sum_{i=1}^nr^i (x^i, \alpha^i)$.
Then $f(x^i) \leq \alpha^i (i=1\cdots, n)$ and
since $f$ is  additive,  it follows that
\begin{eqnarray*}
rf(x)=f(rx)=f(\sum_{i=1}^nr^ix^i)=\sum_{i=1}^nf(r^ix^i)=\sum_{i=1}^n r^if(x^i)\leq \sum_{i=1}^n r^i\alpha^i=r\alpha
\end{eqnarray*}
and so $f(x)\leq \alpha$, which yields that $(x,\alpha)\in  {\rm epi}\,f$. This proves that ${\rm epi}\,f$ is a convex cone.
\cvd

% $a$ is positively homogeneous since  $a(kx)=a(\overbrace{x+\cdots+x}^k)=ka(x)$ for any $x\in X$ and $k\in \mathbb{N}\cup\{0\}$.

\begin{proposition} \label{Pro:additive eq}  {\it Let $X$ be a commutative group and $f:X\rightarrow \mathbb{R}\cup\{\pm \infty\}$. Then the following statements are equivalent:
\begin{itemize}
\item[\rm(i)]  $f$ is  additive.

\item[\rm(ii)] $\pm f$ are subadditive.

\item[\rm(iii)] $f$  is generalised affine such that  either $f:X\rightarrow \mathbb{R}$ with $f(0)=0$, or $f\equiv+\infty$, or $f\equiv -\infty$.
\end{itemize}
Specially, if $f:X\rightarrow \mathbb{R}$ is  additive, then it is $\mathbb{N}$-sublinear with  with $f(0)=0$.
}
\end{proposition}
{\bf Proof.} (i) $\Longleftrightarrow$ (ii)  is clear. It suffices to show (i) $\Longleftrightarrow$ (iii).

(i) $\Longrightarrow$ (iii). It follows from Proposition \ref{Pro:additive char} (i).

(iii) $\Longrightarrow$ (i). If $f\equiv+\infty$, or $f\equiv -\infty$, then (i) holds. Assume that $f:X\to \mathbb{R}$  is affine with $f(0)=0$.
Let $x, y\in X$ and set $z:=x+y$. Then $2z=2(x+y)=(2x)+(2y)$. Since $f$ is affine with $f(0)=0$, from Proposition \ref{Pro:4.4} (ii) we have $f$ is positively homogeneous and taking in account Proposition \ref{Pro:4.2a} it follows that
$
2f(z)=f(2x)+f(2y)=2f(x)+2f(y)
$
and so
$f(z)=f(x)+f(y)$,
yielding that $f$ is additive, i.e., (i) holds.

Specially, if $f:X\rightarrow \mathbb{R}$ is  additive, then from previous equivalence one has $f$ is subadditive with $f(0)=0$. It suffices to show $f$ is positively homogeneous.
Let $x\in X$ and $k\in \mathbb{N}\cup\{0\}$. If $k\in \mathbb{N}$, then from the fact that $f$ is additive we have  $f(kx)=f(\overbrace{x+\cdots+x}^{k})=kf(x)$.
If $k=0$, then from $f:X\rightarrow \mathbb{R}$ and $f(0)=0$ it follows  that
$f(0x)=f(0)=0=0f(x)$. Both cases imply that $f(kx)=kf(x)$.

This completes the proof.
\cvd

\begin{remark} \label{Rm:affineadditive}\it By Proposition \ref{Pro:additive eq} it follows that $f:X\to\mathbb{R}$ is  affine if and only if $f=\psi + c$, where $\psi:X\to\mathbb{R}$ is additive and $c=f(0)\in\mathbb{R}$.  In fact, if $f$ is affine, then setting $\psi:=f-f(0)$ allows that $\psi(0)=0$  and since  the sum of an affine function  and a constant is affine, by Proposition \ref{Pro:additive eq}  we have that $\psi$ is additive. The reverse implication is obvious from the facts that the sum of an affine function and a constant is affine, and an additive function is affine, in view of Proposition \ref{Pro:additive char} (i).

\end{remark}

The following follows immediately from Proposition \ref{Pro:refinementof4.3} and Remark \ref{Rm:affineadditive}.
\begin{proposition} \label{Pro:generalisedaffine char}  {\it Let $X$ be a commutative group and $f:X\rightarrow \mathbb{R} \cup \{\pm\infty\}$ be a generalised affine function. Then either $f:X\to\mathbb{R}$ is  affine with  $f=\psi + f(0)$, where $\psi:X\to\mathbb{R}$ is additive, or $f\equiv+\infty$, or $f\equiv-\infty$.
}
\end{proposition}

\begin{proposition} \label{Pro:exconvfun} {\it Let $A$ be a nonempty  subset of a monoid $X$ and $h:A\to \mathbb{R}$  be a function. Define
the function $g:X\to \mathbb{R}\cup \{+\infty\}$ by
\begin{eqnarray*}
 g(x):=\begin{cases} h(x), \quad {\rm if }\ x\in A, \cr +\infty, \quad  \,{\rm otherwise}. \cr
 \end{cases}
 \end{eqnarray*}
Then the following statements are true:
\begin{itemize}
\item[\rm(i)] Assume that $A$ and $h$ are convex. Then $g$ is convex.

\item[\rm(ii)] Assume that $A$ is a convex cone  with $0\in A$  and $h$ is $\mathbb{N}$-sublinear (resp., subadditive, positively homogeneous). Then $g$ is $\mathbb{N}$-sublinear (resp., subadditive, positively homogeneous).
\end{itemize}
}
\end{proposition}
{\bf Proof.} (i) Let $x\in X,  n, r, r^i\in \mathbb{N}, x^i\in X (i=1,\cdots, n)$ such that
$rx=\sum_{i=1}^nr^i x^i$ and $r=\sum_{i=1}^nr^i$. If $x^i\in A$ for each $i=1,\cdots, n$, then from the facts that $A$ is convex and $h:A\to \mathbb{R}$  is a convex function it follows that $x\in A$ and $rg(x)=rh(x)\leq \sum_{i=1}^nr^i  h(x^i)=\sum_{i=1}^nr^i  g(x^i)$. If $x^i\not\in A$, i.e., $g(x^i)=+\infty$, for some $i=1,\cdots, n$, then $rg(x) \leq  \sum_{i=1}^nr^i  g(x^i)$ holds trivially. This shows that  $g:X\rightarrow \mathbb{R} \cup \{+\infty\}$ is convex.

(ii)  It suffices to show that if $h$ is $\mathbb{N}$-sublinear, then so is $g$.
Since $A$ is a convex cone  with $0\in A$  and $h:A\to \mathbb{R}$ is $\mathbb{N}$-sublinear, we have $h(0)=h(2\cdot0)=2h(0)$ and so $h(0)=0$.
Since $A$ is a convex cone, we first note that if $x,y\in A$ and $n\in \mathbb{N}$, then $x+y, nx\in A$, in view of Propositions \ref{Pro:3.10} and \ref{Pro:3.6} (iii).
Let $x,y\in X$ and  $n\in \mathbb{N}\cup\{0\}$.
If $x, y\in A$, then it is easy to see that $g(x+y)=h(x+y)\leq h(x)+h(y)=g(x)+g(y)$ and $g(nx)=h(nx)=nh(x)=ng(x)$.

Assume that $x\not\in A$ (resp., $y\not\in A$). Then  $g(x+y)\leq +\infty=g(x)+g(y)$ holds trivially. We declare that,
$\mathbf{n}x\not\in A$ (resp., $\mathbf{n}y\not\in A$) for any $\mathbf{n}\in\mathbb{N}$.  If not, i.e., $\mathbf{n}^0x \in A$ (resp., $\mathbf{n}^0y \in A$) for some $\mathbf{n}^0\in\mathbb{N}$,  then from Propositions \ref{Pro:3.10} we have $x \in A$ (resp., $y \in A$), a contradiction. Therefore, $g(nx)=+\infty=ng(x)$ (resp., $g(ny)=+\infty=ng(y)$)
for $n\in\mathbb{N}$. For $n=0$, $g(nx)=g(0)=h(0)=0=0\cdot (+\infty)=n g(x)$ (resp., $g(ny)=g(0)=h(0)=0=0\cdot (+\infty)=n g(y)$).
\cvd

\begin{remark}\label{rem:4.3}{\it In {\rm\cite[Proposition 4.8 (ii)]{LM24}}, we  also assumed that $0\cdot (+\infty) =0$. However,  Borwein and Giladi assumed that $0\cdot (+\infty)=+\infty$ in {\rm\cite[Remark 9]{BG18}},
and they proved {\rm\cite[Theorem 2]{BG18}} without using this assumption (see also, Lemma \ref{LM:5.1}).
}
\end{remark}

\begin{definition}\label{Def:4.6} {\rm (Convex-like and concave-like function)} {\rm  \cite[Definition 4.2]{LM23}} {\it Let $X, Y$ be  semimodules over a semiring $R$. Assume that $Y$ is equipped with a partial order $\leq$. A function $f:X\to Y\cup\{+\infty\}$  is said to be convex-like (resp., concave-like), if for every $n\in \mathbb{N}$, every $r^i\in R\backslash \{0\}$ and every
$x^i\in X (i=1,\cdots,n)$, there exists $z\in X$ such that
$(\sum_{i=1}^nr^i) f(z)\leq \sum_{i=1}^nr^if(x^i)$ (resp., $(\sum_{i=1}^nr^i) f(z)\geq \sum_{i=1}^nr^if(x^i)$).
}
\end{definition}

%{\color{red} Why in the specification of the previous definition do we write {\rm "(Convex-like function on monoid)"} while $X$ is a semimodule?}

\begin{proposition} \label{Pro:4.5} {\rm  \cite[Proposition 4.1]{LM23}}  {\it  Let $X, Y$ be  semimodules over a semiring $R$. Assume that $Y$ is equipped with a compatible partial order $\leq$ and $X$ is divisible. If  $f:X\to Y\cup\{+\infty\}$ is convex (resp., concave), then it is convex-like (resp., concave-like).
}
\end{proposition}

\begin{proposition} \label{Pro:4.6}  {\rm  \cite[Proposition 4.2]{LM23}}  {\it  Let $X$ be a monoid and $Y$ be a commutative group equipped  with a partial order $\leq$ defined by a convex cone $D\subseteq Y$ with $0\in D$. Let $f:X\to Y$.
Then the following statements are true:
\begin{itemize}
\item[\rm(i)] Assume that $Y$ is divisible. If $f(X)+D$ is convex, then  $f$ is convex-like.

\item[\rm(ii)] Assume that $(Y, \leq)$ satisfies that for every $m\in \mathbb{N}, y^1, y^2\in Y$, $ my^1\leq m y^2\Longrightarrow y^1\leq y^2$. If $f$ is convex-like, then $f(X)+D$ is convex.
\end{itemize}
}
\end{proposition}

\section{{New separation theorems between convex sets in groups}}
\setcounter{equation}{0}

In this section, we shall prove several new separation theorems between convex sets in groups by exploiting the following extension of additive functions from a subgroup to the whole group.

\begin{lemma}\label{LM:5.1}{\rm(Extension of additive function)} {\it  Let $X'$ be a convex subgroup of a commutative group $X$, and let $a: X'\to \mathbb{R}$ be additive. Then  there exists an additive function $A: X\to \mathbb{R}$ such that
$A|_{X'}=a$, i.e., $A(x)=a(x)$ for any $x\in X'$. }
\end{lemma}
{\bf Proof.} Define $\mathcal{P}$ the pairs of $(H, h)$ such that $X'\subseteq H\subseteq X$ and $h: H\to \mathbb{R}$ is additive with $h|_{X'}=a$, where $H$ is a convex subgroup of $X$.
Define the following partial order on $\mathcal{P}$ by $(H, h)\leq (H', h')\Longleftrightarrow H\subseteq H'$ and $h'|_H=h$.
Let  $\{(H^i, g^i)\}\subseteq \mathcal{P}$ be a chain. Set $H^0:=\bigcup_i H^i$. Observe that $H^0$ is a convex subgroup of $X$, since the convex subgroups $H^i$ are linearly ordered by inclusion.
Now, define $g^0(h):=g^i(h), h\in H^i$. Since the $H^i$ form a chain, $g^0$ is well defined and so additive. Moreover, $g^0|_{X'}=a$. As a consequence,  $(H^0, g^0)$ is an upper bound for the chain.
Then, this chain  has a maximal element by Zorn's lemma, say, $(H,g)$. We declare that $H=X$. Suppose to the contrary that there is some $x^0\in X\backslash H$. Since $H$ is a convex subgroup with $x^0\not\in H\ni 0$,
from Proposition \ref{Pro:3.3} we have $x^0$ has infinite order and so by Lemma \ref{LM:2.2}, we have  $kx^0=0$ implies that $k=0$, where $k\in \mathbb{Z}$.

Observe that $H$ and $\mathbb{Z}x^0:=\{kx^0: k\in \mathbb{Z}\}$ are subgroups and so is $\mathbf{H}:=H+\mathbb{Z}x^0=\{h+kx^0: h\in H, k\in \mathbb{Z}\}$. Clearly, $X'\subseteq H\subseteq \mathbf{H}$. Let $\mathbf{K}:=H\cap \mathbb{Z}x^0$. Then, $\mathbf{K}$ is a subgroup of $\mathbb{Z}x^0$ and so cyclic.
If $\mathbf{K}=\{0\}$, then choose any $r\in \mathbb{R}$, for instance $r=0$. If $\mathbf{K}\neq\{0\}$, then there is some $d\in \mathbb{N}$ such that $\mathbf{K}=\mathbb{Z}(dx^0)$. In particular, $dx^0\in H$.
Let $e:=g(dx^0)$ and set $r:=\frac{e}{d}$. Now, we define $\bar g: \mathbf{H}\to \mathbb{R}$ by
\begin{eqnarray*}
\bar g(x):=g(h)+kr,\quad x:=h+kx^0\in \mathbf{H}.
\end{eqnarray*}

We now show that $\bar g$ is well defined. Let $x:=h+kx^0, y:=h'+k'x^0\in \mathbf{H}$ and suppose that $x=y$. Then $(k-k')x^0=h'-h\in H$ and so $(k-k')x^0\in H\cap \mathbb{Z}x^0=\mathbf{K}$.
If $\mathbf{K}=\{0\}$, then $(k-k')x^0=0$ and so $h=h'$ and $k=k'$,  by Lemma \ref{LM:2.2}, since $x^0$ has infinite order. This shows that $\bar g(x)=g(h)=g(h')=\bar g(y)$.  If $\mathbf{K}\neq\{0\}$, then $\mathbf{K}=\mathbb{Z}(dx^0)$ for some $d\in \mathbb{N}$. As a consequence, $(k-k')x^0=dsx^0$ for some $s\in \mathbb{Z}$ and so $(k-k'-ds)x^0=0$, which yields that $k-k'=ds$, by Lemma \ref{LM:2.2}, since $x^0$ has infinite order. Therefore,
\begin{eqnarray*}
g(h'-h)=g((k-k')x^0)=g(dsx^0)=sg(dx^0)=sdr=(k-k')r
\end{eqnarray*}
and so
\begin{eqnarray*}
g(h')=g(h)+(k-k')r,
\end{eqnarray*}
which implies that
\begin{eqnarray*}
\bar g(x)=g(h)+kr=g(h')+k'r=g(h')=\bar g(y).
\end{eqnarray*}
This proves that $\bar g$ is well defined.

Next, we show that $\bar g$ is additive.  Let $x:=h+kx^0, y:=h'+k'x^0\in \mathbf{H}$. Then
\begin{eqnarray*}
\bar g(x+y)&=&\bar g((h+h')+(k+k')x^0)\\
&=&g(h+h')+(k+k')r\\
&=&g(h)+g(h')+kr+k'r\\
&=&(g(h)+kr)+(g(h')+k'r)\\
&=&\bar g(x)+\bar g(y),
\end{eqnarray*}
which proves that $\bar g$ is additive. Clearly, $\bar g|_{H}=g$ and $\bar g|_{X'}=a$. Then $(H, h)\leq (\mathbf{H}, \bar g)$ and $H\not \ni x^0\in \mathbf{H}$, a contradiction with the fact that
$(H,g)$ is the  maximal element. Therefore, $H=X$. Now, setting $A:=g$ yields that $A: X\to \mathbb{R}$ is an additive function such that
$A|_{X'}=a$.
\cvd

Let $A$ and $B$ be two nonempty subsets of a monoid $X$.
Recall that $A$ and $B$ are separable
if there is a nonzero  additive function $a:X\to \mathbb{R}$ such that
$\sup_{x\in A}a(x) \leq \inf_{x\in B}a(x)$, i.e., $a(x) \leq  a(y)$ for any $x\in A$ and $y\in B$.
If, in addition,
$\inf_{x\in A}a(x) < \sup_{x\in B}a(x)$,
i.e. there are $x^1\in A$  and $x^2\in B$ such that $a(x^1)<a(x^2)$, then $A$ and $B$ are properly separable.

\begin{lemma} \label{Thm:Separation} { \rm \cite[Theorem 5.4]{LM24} \cite[Theorem 6.1]{LM25}}  {{\it Let $C$ be a convex subset of a   commutative group $Z$} with   ${\rm qcore}\, C\not=\emptyset$
and let $x^0\not\in C$. Then $C$ and $x^0$ are properly separable.
}
\end{lemma}

%{\color{blue}
%\begin{proposition}\label{Pro:5.1}  {\it Let $X$ be a commutative group and $x^0\in X$.  Then $X^0:=\{kx^0: k\in \mathbb{Z}\}$ is a convex subgroup of $X$.
%}
%\end{proposition}
%{\bf Proof.}\, It is clear that $X^0$ is a subgroup of $X$. It sufficed to show $X^0$ is convex. To this aim, let $x\in X$, $n, r, r^1, \cdots, r^n \in \mathbb{N}$ and
%$k^1x^0,\cdots, k^nx^0 \in X^0 $ such that
%\begin{eqnarray*}
%rx=\sum_{i=1}^nr^i k^ix^0,\,\, r=\sum_{i=1}^nr^i.
%\end{eqnarray*}
%}

\begin{proposition} \label{Pro:5.1} {\it Let $C, D$ be subsets of a monoid $X$. Then the following statements are equivalent:
\begin{itemize}
\item[\rm(i)] $C$ and $D$ are (resp., properly) separable.

\item[\rm(ii)] ${\rm conv}\,C$ and ${\rm conv}\, D$ are (resp., properly) separable.

\item[\rm(iii)] $C$ and ${\rm conv}\,D$ are (resp., properly) separable.

\item[\rm(iv)] ${\rm conv}\,C$ and $D$ are (resp., properly) separable.
\end{itemize}
}
\end{proposition}
{\bf Proof.}\, It suffices to show (i) $\Longleftrightarrow$ (ii). Similarly, we can show (ii) $\Longleftrightarrow$ (iii) and (iii) $\Longleftrightarrow$ (iv).
Observe that
\begin{eqnarray*}
C\, {\rm and}\ D\, {\rm are\, properly\, separable}&\Longleftrightarrow & {\rm  there\, is\, a\, nonzero\,  additive\, function}\,a:X\to \mathbb{R}\, {\rm such\, that}\\
&&\sup_{x\in C}a(x) \leq \inf_{x\in D}a(x)\,{\rm and}\,\inf_{x\in C}a(x) < \sup_{x\in D}a(x)\\
&\overset{\rm (By\, Proposition\, \ref{Pro:4.4+})}{\Longleftrightarrow}& {\rm  there\, is\, a\, nonzero\,  additive\, function}\,a:X\to \mathbb{R}\, {\rm such\, that}\\
&&\sup_{x\in {\rm conv} C}a(x) \leq \inf_{x\in {\rm conv} D}a(x)\,{\rm and}\,\inf_{x\in {\rm conv} C}a(x) < \sup_{x\in {\rm conv}D}a(x)\\
&\Longleftrightarrow & {\rm conv}\,C\, {\rm and}\ {\rm conv}\,D\, {\rm are\, properly\, separable}.
\end{eqnarray*}
Similarly, we can show that $C$ and $D$ are separable if and only if ${\rm conv}\,C$ and ${\rm conv}\, D$ are separable.
\cvd

From Proposition \ref{Pro:5.1} one has the following.
\begin{corollary}   \label{Cor:5.1-}  {\it Let $C$ be a subset of a monoid $X$ and $x^0\in X$. Then the following statements are equivalent:
\begin{itemize}
\item[\rm(i)] $C$ and $x^0$ are (resp., properly) separable.

\item[\rm(ii)] ${\rm conv}\,C$ and ${\rm conv}\, \{x^0\}$ are (resp., properly) separable.

\item[\rm(iii)] $C$ and ${\rm conv}\,\{x^0\}$ are (resp., properly) separable.

\item[\rm(iv)] ${\rm conv}\,C$ and $x^0$ are (resp., properly) separable.
\end{itemize}
}
\end{corollary}
%{\bf Proof.}\, From Proposition \ref{Pro:singletonconv} we have ${\rm conv}\, \{x^0\}=\{x^0\}$. Therefore, the conclusion follows immediately from Proposition \ref{Pro:5.1}.
%\cvd

\begin{theorem} \label{Thm:5.2} {\rm (Proper separation between a point and a convex set I)}  {\it Let $C$ be a convex subset of a   commutative group $X$ with   ${\rm qr core}\, C\not=\emptyset$
and let $X\ni x^0\not\in C$. Then $C$ and $x^0$ are properly separable.
}
\end{theorem}
{\bf Proof.}\, Since ${\rm qr core}\, C\not=\emptyset$, let $y^0\in {\rm qr core}\, C$ and set $X':={\rm cone}\,(C-y^0)$. Then $C-y^0\subseteq X'$ and from  Definition \ref{Def:3.5} (iii), $X'$ is a convex  subgroup of $X$. We prove the theorem in the following two cases.

Case 1. Assume that $x^0-y^0\not\in X'$. Then Since $X'$ is a convex  subgroup of $X$, $0\in X'$ and from Proposition \ref{Pro:3.3} we have $x^0-y^0$ has infinite order. Observe that $X'$ and $\mathbb{Z}(x^0-y^0):=\{k(x^0-y^0): k\in \mathbb{Z}\}$ are subgroups of $X$ and so is $X'+\mathbb{Z}(x^0-y^0)=\{x'+k(x^0-y^0): x\in X', k\in \mathbb{Z}\}$. Let $\mathbf{X}:={\rm cone}\,(X'+\mathbb{Z}(x^0-y^0))$. Then $X'\subseteq \mathbf{X}, \mathbb{Z}(x^0-y^0)\subseteq \mathbf{X}$ and from Proposition \ref{Pro:conesubgroup}, we have $\mathbf{X}$ is a convex subgroup of $X$.
Now, define $a: \mathbf{X}\to \mathbb{R}$ by
\begin{eqnarray*}
a(x):= \frac{\sum_{i=1}^{n} r^ik^i}{r},\quad x \in \mathbf{X},
\end{eqnarray*}
where $x$ satisfies that $rx=\sum_{i=1}^n r^i(x^i+k^i(x^0-y^0))=\sum_{i=1}^n r^ix^i+\sum_{i=1}^n r^ik^i(x^0-y^0)$
for $n,r, r^i\in \mathbb{N}$, $k^i\in \mathbb{Z}$ and $x^i\in X' (i=1,\cdots,n)$.

We first show that the function $a: \mathbf{X}\to \mathbb{R}$ is well defined. To this aim, assume that $x\in \mathbf{X}$ satisfies that
$rx=\sum_{i=1}^{n'} r^i(x^i+k^i(x^0-y^0))=\sum_{i=1}^{n'} r^ix^i+\sum_{i=1}^{n'}r^i k^i(x^0-y^0)$ for $n',r, r^i\in \mathbb{N}$, $k^i\in \mathbb{Z}$, $x^i\in X'(i=1,\cdots,n')$, and
$sx=\sum_{j=1}^{m'} s^j(y^j+t^j(x^0-y^0))=\sum_{j=1}^{m'} s^jy^j+\sum_{j=1}^{m'}s^j t^j(x^0-y^0)$
for $m',s,s^j\in \mathbb{N}$, $t^j\in \mathbb{Z}$ and $y^j\in X'  (j=1,\cdots,m')$.
For the first case, we have $a(x)=\frac{\sum_{i=1}^{n'} r^ik^i}{r}$, and for the second one, we have  $h(x)=\frac{\sum_{j=1}^{m'} s^jt^j}{s}$. It suffices to show
$h(x)=\frac{\sum_{i=1}^{n'} r^ik^i}{r}=\frac{\sum_{j=1}^{m'} s^jt^j}{s}$.
In fact, from the facts that $rx=\sum_{i=1}^{n'} r^ix^i+\sum_{i=1}^{n'}r^i k^i(x^0-y^0))$ and
$sx=\sum_{j=1}^{m'} s^jy^j+\sum_{j=1}^{m'}s^j t^j(x^0-y^0)$ it follows that $srx=\sum_{i=1}^{n'} sr^ix^i+\sum_{i=1}^{n'}sr^i k^i(x^0-y^0)$ and
$rsx=\sum_{j=1}^{m'} rs^jy^j+\sum_{j=1}^{m'}rs^j t^j(x^0-y^0)$. Since $X'$ is a convex  subgroup of $X$, it follows that
\begin{eqnarray*}
(\sum_{i=1}^{n'}sr^i k^i-\sum_{j=1}^{m'}rs^j t^j)(x^0-y^0)=\sum_{i=1}^{n'} sr^ix^i-\sum_{j=1}^{m'} rs^jy^j\in X'
\end{eqnarray*}
and so
\begin{eqnarray*}
\pm(\sum_{i=1}^{n'}sr^i k^i-\sum_{j=1}^{m'}rs^j t^j)(x^0-y^0)\in X'.
\end{eqnarray*}
We declare that $\sum_{i=1}^{n'}sr^i k^i=\sum_{j=1}^{m'}rs^j t^j$. In fact, if $\sum_{i=1}^{n'}sr^i k^i\neq \sum_{j=1}^{m'}rs^j t^j$, then it follows that
$x^0-y^0 \in X'$, in view of Proposition \ref{Pro:3.10}, a contradiction. Thus $\sum_{i=1}^{n'}sr^i k^i=\sum_{j=1}^{m'}rs^j t^j$ and consequently, $a(x)=\frac{\sum_{i=1}^{n'} r^ik^i}{r}=\frac{\sum_{j=1}^{m'} s^jt^j}{s}$.

We next show that the function $a: \mathbf{X}\to \mathbb{R}$ is additive. Let $x\in \mathbf{X}$ be such that
$rx=\sum_{i=1}^{n'} r^i(x^i+k^i(x^0-y^0))=\sum_{i=1}^{n'} r^ix^i+\sum_{i=1}^{n'}r^i k^i(x^0-y^0)$ for $n',r, r^i\in \mathbb{N}$, $k^i\in \mathbb{Z}$, $x^i\in X'(i=1,\cdots,n')$, and let $y\in \mathbf{X}$  be such that
$sy=\sum_{j=1}^{m'} s^j(y^j+t^j(x^0-y^0))=\sum_{j=1}^{m'} s^jy^j+\sum_{j=1}^{m'}s^j t^j(x^0-y^0)$
for $m',s,s^j\in \mathbb{N}$, $t^j\in \mathbb{Z}$ and $y^j\in X'  (j=1,\cdots,m')$.
Then $h(x)=\frac{\sum_{i=1}^{n'} r^ik^i}{r}$ and $h(y)=\frac{\sum_{j=1}^{m'} s^jt^j}{s}$. Since $\mathbf{X}$ is a convex subgroup we have  $x+y\in X'$.
Moreover,
$srx=\sum_{i=1}^{n'} sr^ix^i+\sum_{i=1}^{n'}sr^i k^i(x^0-y^0)$ and $rsy=\sum_{j=1}^{m'} rs^jy^j+\sum_{j=1}^{m'}rs^j t^j(x^0-y^0)$ and so
\begin{eqnarray*}
sr(x+y)=(\sum_{i=1}^{n'} sr^ix^i+\sum_{j=1}^{m'} rs^jy^j)+(\sum_{i=1}^{n'}sr^i k^i+\sum_{j=1}^{m'}rs^j t^j)(x^0-y^0),
\end{eqnarray*}
which yields that
\begin{eqnarray*}
h(x+y)=\frac{\sum_{i=1}^{n'} s r^i k^i+\sum_{j=1}^{m'} r s^j t^j}{sr}=\frac{\sum_{i=1}^{n'} r^i k^i}{r}+\frac{\sum_{j=1}^{m'} s^j t^j}{s}=h(x)+h(y),
\end{eqnarray*}
which implies $h$ is additive. It is easy to check that $a(u)=0$ for any $u\in X'$ and $a(x^0-y^0)=1$ and so $a(x^0)=a(y^0)+1$.

Since $\mathbf{X}$ is a convex subgroup of $X$ and $a: \mathbf{X}\to \mathbb{R}$ is nonzero additive, from Lemma \ref{LM:5.1} there is a nonzero additive function $A: X\to \mathbb{R}$ such that
$A|_{\mathbf{X}}=a$. Observe that $X'\subseteq \mathbf{X}, \mathbb{Z}(x^0-y^0)\subseteq \mathbf{X}$. Then, we have $x-y^0\in X'$ for any $x\in C$ and so $A(x-y^0)=a(x-y^0)=0$, i.e.,  $A(x)=A(y^0)$ for any $x\in C$. Also,
$A(x^0-y^0)=a(x^0-y^0)=1$ and so $A(x^0)=A(y^0)+1$.
As a consequence,
\begin{eqnarray*}
A(x)=A(y^0)<A(y^0)+1=A(x^0),\quad \forall x\in C,
\end{eqnarray*}
which yields that $C$ and $x^0$ are properly separable.

Case 2. Assume that $x^0-y^0 \in X'$. Observe that $X'={\rm cone}\,(C-y^0)$ and $X'$ is a convex  subgroup of $X$.
Now, work in the group $X'$ with $C':=C-y^0$, we have  $X'={\rm cone}\,(C'-0)$ and so $0\in{\rm qcore}_{X'}\,C'$, where ${\rm qcore}_{X'}$ denotes the quasi core with respect to $X'$.
Clearly, $x^0-y^0\not\in C'$ and ${\rm qcore}_{X'}\,C'\neq\emptyset$.
Now, applying Lemma \ref{Thm:Separation} yields that $C'$ and $x^0-y^0$ are properly separable by a nonzero additive function $a: X'\to \mathbb{R}$, that is, $a(x-y^0) \leq a(x^0-y^0)$  for any $ x\in C$, and $l(v^0-y^0)< l(x^0-y^0)$ for some $v^0\in C$. As a consequence, $a(x)-a(y^0)=a(x-y^0) \leq a(x^0-y^0)=a(x^0)-a(y^0)$, i.e.,  $a(x) \leq a(x^0)$ for any $ x\in C$,
and $a(v^0)-a(y^0)=a(v^0-y^0)< a(x^0-y^0)=a(x^0)-a(y^0)$, i.e., $a(v^0)<a(x^0)$ for some $v^0\in C$. This implies that $C$ and $x^0$ are properly separable by a nonzero additive function $a: X'\to \mathbb{R}$.
Now, from Lemma \ref{LM:5.1} there is a  nonzero additive function $A: X\to \mathbb{R}$ such that
$A|_{X'}=a$. Therefore, $C$ and $x^0$ are properly separable by the nonzero additive function $A: X\to \mathbb{R}$. This completes the proof.
\cvd

%\begin{remark}   \label{Rm:5.1}{\it Observe that in Theorem {\rm\ref{Thm:5.2}}, since $\{kx^0: k\in \mathbb{Z}\}$ is a subgroup of $X$, from Proposition  \ref{Pro:conesubgroup} it follows that $X'$ is a convex subgroup of $X$. %We also  have {\color{red}$l:X\to \mathbb{R}$ with $l(0)=0$ and $l(x_0)\geq 1$}  and $\rho_C(x) < +\infty$ for any $x\in {\rm cone}\,C$, in view of Proposition {\rm\ref{Pro:4.2} (ii)}.
%Moreover, if $x^0\not\in {\rm cone}\,C$, then from Remark {\rm\ref{Rm:Mink}} one has $\rho_C(x^0)=+\infty$. As a consequence, in order to assure that $\rho_C(x^0)<+\infty$, we may assume that
%$x^0 \in {\rm cone}\,C$.
%}
%\end{remark}

The following example proposed by Professor C. Ritzenthaler is given to illustrate Theorem \ref{Thm:5.2}. This implies that a convex subset  and a point in a commutative group are not properly separable by a nonzero additive function  without any other additional assumptions, which answers affirmatively the open question proposed by Li and Mastroeni in \cite{LM25}.

\begin{example}\label{Ex:SBPS}{\rm Let $X:=\mathbb{Z}^{\mathbb{N}}$, the group of finitely supported integer sequences (that is, any element $x\in X$ has only finitely nonzero integer coordinates), with unit elements $e^i:=(0,\cdots,0,\stackrel{\small(i)}{1},0,\cdots), i=1, 2,\cdots$. Let $C:=\{0\}\cup \{x!=0:{\rm the\ last\ nonzero\ coordinate\ of}\ x\in X \ {\rm is\ positive}\}$ (the lexicographic cone), and $x^0:=-e^1$.

It is easy to check that $0\in C$, $x^0 \notin C$ and $C$ is a convex cone but not a subgroup in $X$. We next show that ${\rm qr core}\, C=\emptyset$. In fact, $0\not\in {\rm qr core}\, C$, by Proposition \ref{Pro:3.6} (i), since ${\rm cone}\, (C-0)=C$. Let $u:=(u_1, u_2,\cdots)\in C\backslash\{0\}$ and denote by $L(u)$ the coordinate number of the last nonzero coordinate. Then $L(u)>0$ and define $v\in X$ by
\begin{eqnarray*}
 v_i:=\begin{cases} u_i, \quad\  i=1,\cdots,L(u), \cr 1, \quad\,\,\,  \,i=L(u)+1, \cr 0, \quad\,\,\,\,\, i=L(u)+2,\cdots.
 \end{cases}
 \end{eqnarray*}
It is clear that $v\in {\rm cone}\, (C-u)$, but $-v\not\in {\rm cone}\, (C-u)$. This shows that ${\rm cone}\, (C-u)$ is not a subgroup of $X$ for any $C\backslash\{0\}$. This proves that ${\rm qr core}\, C=\emptyset$.

Let $l: X \to \mathbb{R}$ be additive with $l(y) \leq  l(x^0)$ for all $y \in C$.
For any $y \in C$ and $k \in \mathbb{N}$, we have $ky \in C$ and it follows that $kl(y)=l(ky) \leq l(x^0)$ and so $l(y)\leq \frac{1}{k}l(x^0)\to 0$ as $k\to +\infty $, i.e., $l \leq 0$ on $C$. Since $e^n\in C$, $l(e^n) \leq 0$ for every $n \in \mathbb{N}$. Moreover,  $e^{n+1}+t e^n \in C$ for every $t \in \mathbb{Z}$, so $l(e^{n+1})+tl(e^n)=l(e^{n+1}+te^n)\leq 0$ for all $t \in \mathbb{Z}$. Hence $l(e^{n+1})\leq kl(e^n)\leq0$ for every $n, k\in \mathbb{N}$. We declare that $l(e^n)\geq0$ for every $n \in \mathbb{N}$. If not, there is $n'\in N$ such that
$l(e^{n'})<0$, then $kl(e^{n'})\to -\infty$ as $k\to +\infty$, a contradiction with $+\infty<l(e^{n'+1})\leq kl(e^{n'})\to -\infty$. Then $l(e^n)\geq0$ and so $l(e^n)=0$ for every $n \in \mathbb{N}$. Observe that
$x=\sum_{i=1}^{+\infty}e^ix_i$ for every $x\in X$. As a consequence, we have $l=0$ on $X$. This shows that $C$ and $x^0$ are not properly separable by a nonzero additive function.

}
\end{example}

%\begin{comment}
%\begin{theorem} \label{Thm:5.3} {\rm (Proper separation between a point and a convex set II)}  {\it Let $C$ be a nonempty convex subset of a commutative group $X$ with   ${\rm qr core}\, C\not=\emptyset$
%and let $X\ni x^0\not\in C$. Then $C$ and $x^0$ are properly separable.
%}
%\end{theorem}
%{\bf Proof.}\, Let $v^0\in C$. Then $0\in C-v^0$. Since $C$ is convex, from Proposition  \ref{Pro:3.4}, so is $C-v^0$. Observe that $x^0-v^0 \not\in C-v^0$. If not,  $x^0-v^0 \in C-v^0$, i.e., $x^0\in C$, a contradiction.
%Since ${\rm qr core}\, C\not=\emptyset$, from Proposition  \ref{Pro:eqqqrcore} we have ${\rm qr core}\, (C-v^0)={\rm qr core}\, C-v^0\neq \emptyset$.
%Now from Theorem \ref{Thm:5.2}, $C-v^0$ and $x^0-v^0$ are properly separable. Then there is a nonzero additive function $l:X\to \mathbb{R}$  such that
%$ l(y-v^0) \leq l(x^0-v^0)$   for any $y\in C$ and
%$l(0)< l(x^0-v^0)$.
%Since $l:X\to \mathbb{R}$ is additive, it follows that
%\begin{eqnarray*}
%l(y)&=&l(y)-l(v^0)+ l(v^0) =l(y-v^0)+ l(v^0) \\
%&\leq& l(x^0-v^0)+l(v^0)=l(x^0)-l(v^0)+l(v^0)=l(x^0), \forall y\in C,
% \end{eqnarray*}
%and
%\begin{eqnarray*}
%l(v^0)=l(0)+l(v^0)< l(x^0-v^0)+l(v^0)=l(x^0)-l(v^0)+l(v^0)=l(x^0).
% \end{eqnarray*}
%This yields that $C$ and $x^0$ are properly separable.
%This completes the proof.
%\cvd
%\end{comment}

\begin{lemma} \label{LM:6.1} {\rm\cite[Lemma 6.5]{LM25}}  {\it Let $Z$ be a monoid, $a:X\to \mathbb{R}$  be a nonzero additive function and $A\subseteq Z$. Then the following statements are true.
\begin{itemize}
\item[\rm(i)] $a(x)\leq 0$ for any $x\in A$ if and only if $a(x)\leq 0$ for any $x\in {\rm cone}\,A$, or equivalently, $\sup_{x\in A}a(x)\leq 0$ if and only if $\sup_{x\in {\rm cone}\,A}a(x)\leq 0$.

\item[\rm(ii)] $a(x)< 0$ for some $x\in A$ if and only if $a(x)<  0$ for some $x\in {\rm cone}\,A$, or equivalently, $\inf_{x\in A}a(x)<0$ if and only if $\inf_{x\in {\rm cone}\,A}a(x)<0$.
\end{itemize}
}
\end{lemma}

%\begin{definition} \label{Deg:qusircore}  {\it Let $C$ be a nonempty subset of a group $Z$.
%\begin{itemize}
%\item[\rm(i)] {\rm (Quasi core)} \rm \cite[Definition 4.21]{LM24}  The quasi core of  $C$ is defined by
%$
%{\rm qcore}\,C:=\{x\in C:{\rm cone}\,(C-x)=Z\}.
%$

%\item[\rm(ii)] {\rm (Pseudo relative core)}   {\rm \cite[Definition 4.1]{LM25}} The pseudo relative core of  $C$ is defined by
%\begin{eqnarray*}
%{\rm prcore}\,C&:=&\{x\in C:{\rm cone}\,(C-x)\, {\rm is\,\, a\,\, subgroup\,\, of}\, Z\}\\
%&:=&\{x\in C:{\rm cone}\,(C-x)\, {\rm is\,\, a\,\, convex\,\,subgroup\,\, of}\, Z\}.
% \end{eqnarray*}
%}
%\end{itemize}
%}
%\end{definition}
%Note that the second equality in Definition \ref{Deg:qusircore} (ii) follows from Proposition \ref{Pro:3.6} (i).
%Clearly, ${\rm qcore}\,C\subseteq {\rm prcore}\,C\subseteq C$.

We have the following interesting proper separation theorem between a point and a set in commutative groups, in which  we do not need to assume $C$ is convex.
% and ${\rm prcore}\, C\ne\emptyset $.

\begin{theorem} \label{Thm:5.3+}{\rm (Proper separation between a point and a convex set II)} {\it Let $C$ be a nonempty subset of a commutative group, $x^0\in C$ and $x^0\not\in{\rm prcore}\,C\neq \emptyset$. Then $C$ and $x^0$ are properly separable.}
\end{theorem}

\noindent {\bf Proof.}\, Since $x^0\in C$ and $x^0\not\in{\rm prcore}\,C$, ${\rm cone}\, (C-x^0)$ is not a subgroup of $X$ and so there exists $\hat x\in {\rm cone}\,(C-x^0)$ such that $-\hat x\not\in {\rm cone}\,(C-x^0)$. Since ${\rm prcore}\,C\neq \emptyset$, let $y^0\in {\rm prcore}\,C$. Then ${\rm cone}\, (C-y^0)$ is a subgroup of $X$, $x^0-y^0\in C-y^0\subseteq {\rm cone}\, (C-y^0)$ and so $y^0-x^0=-(x^0-y^0)\in {\rm cone}\, (C-y^0)$, which implies that  $C-x^0=(y^0-x^0)+(C-y^0) \subseteq {\rm cone}\, (C-y^0)$ and thus ${\rm cone}\, (C-x^0)\subseteq {\rm cone}\, (C-y^0)$. On the other hand, we have $C-y^0\subseteq {\rm cone}\, (C-x^0)-(y^0-x^0)\subseteq {\rm cone}\, (C-y^0)+{\rm cone}\, (C-y^0)={\rm cone}\, (C-y^0)$, where the equality follows from Proposition \ref{Pro:3.6} (iii). Now, from Proposition \ref{Pro:3.6} (i) it follows that
${\rm cone}\, (C-y^0)\subseteq {\rm cone}\,({\rm cone}\, (C-x^0)-(y^0-x^0))\subseteq {\rm cone}\, (C-y^0)$ and so ${\rm cone}\,({\rm cone}\, (C-x^0)-(y^0-x^0))={\rm cone}\, (C-y^0)$, which yields that
$y^0-x^0\in {\rm prcore}\,({\rm cone}\, (C-x^0))$, i.e., ${\rm prcore}\,({\rm cone}\, (C-x^0))\neq\emptyset$.
Now from Theorem \ref{Thm:5.2},  ${\rm cone}\,(C-x^0)$ and $-\hat x$ are properly separable, i.e.,  there exists a nonzero additive function $l:X\to \mathbb{R}$ such that
\begin{equation}\label{la:a1}
 l(y) \leq l(-\hat x)=-l(\hat x), \,\,  \forall y\in  {\rm cone}\,(C-x^0)
\end{equation}
and
\begin{equation}\label{la:a2}
 l(\hat y) < l(-\hat x)=-l(\hat x)
\end{equation}
for some $\hat y\in {\rm cone}\,(C-x^0)$.

Note that, since  $0\in  {\rm cone}\,(C-x^0)$,  \eqref{la:a1} implies $0=l(0)\le-l(\hat x)$ and so $l(\hat x)\le 0$. We consider the following two cases $l(\hat x)=0$ and $l(\hat x)<0$.

(i) Assume $l(\hat x)=0$. By \eqref{la:a1}, we obtain that
\begin{equation*}
l(x)-l(x^0)=l(x-x^0) \leq 0 , \,\,  \forall x\in C
\end{equation*}
and from the fact that $l$ is additive we have
\begin{equation}\label{la:a3}
l(x)\leq l(x^0), \quad  \forall x\in C.
\end{equation}
Moreover,  by Lemma  \ref{LM:6.1} (ii), from \eqref{la:a2} it follows that there exists $\tilde x\in C$ such that $l(\tilde x-x^0) < 0$, i.e., $l(\tilde x)< l(x^0)$, since $l$ is additive. This, along with
\eqref{la:a3}, yields that $C$ and $x^0$ are properly separable.

(ii) Assume $l(\hat x)<0$.
% and strict inequality holds  in \eqref{la:a1} for one $\hat y\in {\rm cone}\,(C-x^0)$.
Note that $l$ is additive.
Then, by \eqref{la:a1} and Proposition \ref{Pro:3.6} (iii), we have
\begin{equation*}
r l(x-x^0)=l(r(x-x^0)) \le -l(\hat x),\,\,\forall r\in\mathbb{N}, \ \forall x\in C
\end{equation*}
and so
\begin{equation*}
l(x)-l(x^0)=l(x-x^0)\le \frac{-l(\hat x)}{r}\to 0,\,\, {\rm as}\,\, r\to +\infty,\,\, \forall x\in C,
\end{equation*}
which yields that \eqref{la:a3} holds. Since $l(\hat x)<0$ and $\hat x\in {\rm cone}\,(C-x^0)$, by Lemma  \ref{LM:6.1} (ii), there exists $\tilde x\in C$ such that $l(\tilde x-x^0) < 0$, i.e., $l(\tilde x)< l(x^0)$, since $l$ is additive.  This, along with \eqref{la:a3}, yields that $x^0$ and $C$ are properly separable. \cvd

\begin{theorem} \label{Thm:5.3++}{\rm (Proper separation between a point and a convex set III)} {\it Let $C$ be a nonempty convex subset of a commutative group $X$ and  $X\ni x^0\not\in{\rm prcore}\,C\neq\emptyset$. Then $C$ and $x^0$ are properly separable.}
\end{theorem}
{\bf Proof.}\, If $x^0\not\in C$, then the conclusion follows immediately from Theorem \ref{Thm:5.2}. If  $x^0 \in C$, then the conclusion also follows immediately from Theorem \ref{Thm:5.3+}.
\cvd

\begin{remark}   \label{Rm:5.1+}{\rm  If ${\rm qcore}\,C\neq \emptyset$, then from  {\rm \cite[Propositiion 4.3]{LM25}} one has
${\rm qcore}\,C={\rm prcore}\,C$. So, if we assume that ${\rm qcore}\,C\neq \emptyset$, then  Theorem \ref{Thm:5.3++} collapses to {\rm \cite[Theorem 6.2]{LM25}}, which is stronger than {\rm \cite[Theorem 5.4]{LM24}}, where
we assumed that $x^0\not\in C$ instead of $x^0\not\in{\rm qcore}\,C$.
}
\end{remark}

We next show that the condition $x^0\not\in {\rm prcore}\,C $ is also necessary for the proper separation between $C$ and $x^0$.

\begin{proposition}  \label{Pro:5.1+0}{\it Let $C$ be a convex subset of a commutative group $X$  and $x^0\in X$. Assume that $C$ and $x^0$ are properly separable. Then  $x^0\not\in {\rm prcore}\,C$.
}
\end{proposition}
{\bf Proof.}\, If $x^0\not\in\, C$ or ${\rm prcore}\,C=\emptyset$, then the conclusion is obvious. Assume that $x^0\in C$ and ${\rm prcore}\,C\neq\emptyset$. By contradiction, assume that $x^0\in {\rm prcore}\,C$, i.e., ${\rm cone}\,(C-x^0)$ is a subgroup of $X$.
Since $C$ and $x^0$ are properly separable, there exist a nonzero additive function $l:X\to \mathbb{R}$ and $ \alpha\in\mathbb{R}$ such that
\begin{equation}\label{S1}
l(x^0)\ge\alpha,\qquad l(x)\le\alpha,\quad\forall x\in C
\end{equation}
and moreover, there exists $\bar x\in C$ such that $l(\bar x) < \alpha$, which implies
\begin{equation}\label{S2}
l(\bar x-x^0)=l(\bar x) -l(x^0) < \alpha -\alpha =0.
\end{equation}
Setting $z:=x-x^0\in C-x^0$, from \eqref{S1} it follows that
\begin{eqnarray*}
l(z)=l(x) -l(x^0)\le 0,\quad\forall z\in C-x^0,
\end{eqnarray*}
and applying Lemma \ref{LM:6.1} (i) yields that
\begin{eqnarray*}
l(z)\le 0,\quad\forall z\in {\rm cone}\,(C-x^0).
\end{eqnarray*}
Since  ${\rm cone}\,(C-x^0)$ is a subgroup of $X$, $z\in {\rm cone}\,(C-x^0)$ implies that $-z\in {\rm cone}\,(C-x^0)$ and by the previous inequality, we get
\begin{eqnarray*}
l(z) = 0,\quad\forall z\in {\rm cone}\,(C-x^0).
\end{eqnarray*}
In particular, $l(z) = 0$ for any $z\in C-x^0$, which contradicts \eqref{S2}. This completes the proof.
\cvd

\begin{theorem} \label{Thm:5.4-} {\rm (Proper separation between  convex sets I)}  {\it Let $C$ and $D$ be nonempty convex subsets of a commutative group $X$  such that $C-D$ is convex with  $C\cap D = \emptyset$ and ${\rm prcore}\,(C-D)\neq\emptyset$.
 Then  $C$ and $D$ are properly separable.
}
\end{theorem}
{\bf Proof.}\, Since $C\cap D = \emptyset$, one has $0\not\in C-D$. Now, Theorem \ref{Thm:5.2} yields that $C-D$ and $0$ are properly separable. Then there is a nonzero additive function  $l:X\to \mathbb{R}$ such that
$l(y) \leq l(0)$  for any $y \in C-D$
and
$
l(\bar y) < l(0)
$
for some $\bar y\in C-D$.
Then
\begin{eqnarray}\label{5.1}
 l(u-v) \leq l(0), \forall u \in C,\, \forall v\in D,
 \end{eqnarray}
and
\begin{eqnarray}\label{5.2}
 l(\bar u-\bar v) < l(0)
 \end{eqnarray}
for some $\bar u\in C$ and $\bar v\in D$.
Since $l:X\to \mathbb{R}$ is additive, from  \eqref{5.1} and \eqref{5.2} we have
 \begin{eqnarray*}
 l(u)=l(u)-l(v)+l(v) =l(u-v)+l(v) \leq l(0)+l(v)=l(v), \forall u \in C,\, \forall v\in D,
 \end{eqnarray*}
and
\begin{eqnarray*}
l(\bar u)= l(\bar u)-l(\bar v)+l(\bar v) =l(\bar u-\bar v)+l(\bar v) < l(0)+l(\bar v)=l(\bar v)
 \end{eqnarray*}
for some $\bar u\in C$ and $\bar v\in D$.
This proves that $C$ and $D$ are properly separable.
\cvd

\begin{theorem} \label{Thm:5.4} {\rm (Proper separation between  convex sets II)}  {\it Let $C$ and $D$ be nonempty convex subsets of a uniquely divisible commutative group $X$ with  $C\cap D = \emptyset$ and ${\rm prcore}\,(C-D)\neq\emptyset$. Then  $C$ and $D$ are properly separable.
}
\end{theorem}
{\bf Proof.}\, Since $X$ is uniquely divisible and $C$ and $D$ are convex, from Proposition \ref{Pro:3.4+} one has  $C-D$ is  convex.
From Theorem \ref{Thm:5.4-} it follows that $C$ and $D$ are properly separable.
\cvd

\begin{corollary} \label{Thm:5.5} {\rm (Proper separation between   sets)}  {\it Let $C$ and $D$ be nonempty subsets of a uniquely divisible commutative group $X$ with ${\rm conv}\,C\cap {\rm conv}\,D = \emptyset$ and ${\rm prcore}\,({\rm conv}\,C-{\rm conv}\,D)\neq\emptyset$.
 Then  $C$ and $D$ are properly separable.}
\end{corollary}
{\bf Proof.}\, From Theorem \ref{Thm:5.4} we have ${\rm conv}\,C$ and ${\rm conv}\,D$ are properly separable. Now Proposition \ref{Pro:5.1} allows that $C$ and $D$ are properly separable.
\cvd

\vskip0.2cm

In the following, we provide proper separation between  convex cones in groups, in which we remove the assumption that the divisible commutative group $X$ is unique, compared with proper separation between  convex sets in groups in Theorem \ref{Thm:5.4}.

\begin{theorem} \label{Thm:5.6}  {\rm (Proper separation between  convex cones I)}  {\it Let $C$ and $D$ be nonempty convex cones of a  divisible commutative group $X$ with $C\cap D = \emptyset$ and ${\rm prcore}\,(C-D)\neq\emptyset$.
Then there is a nonzero additive function  $l:X\to \mathbb{R}$  such that
\begin{eqnarray}\label{5.4}
l(u) \leq 0 \leq l(v), \forall u\in C,\, \forall v\in D,
 \end{eqnarray}
or equivalently, $l\in (-C^*)\cap D^*$,
and
\begin{eqnarray}\label{5.5}
l(\bar u) <  l(\bar v)
 \end{eqnarray}
for some $\bar u\in C$ and $\bar v\in D$, where either $l(\bar u)<0$ or $l(\bar v)>0$.
}
\end{theorem}
{\bf Proof.}\, Since $X$ is divisible and $C$ and $D$ are convex, from Proposition \ref{Pro:3.9} we have $C-D$ is a convex cone.
From Theorem \ref{Thm:5.4-}, there is a nonzero additive function  $l:X\to \mathbb{R}$ such that
\begin{eqnarray}\label{5.6}
l(u) \leq l(v), \forall u\in C,\, \forall v\in D,
 \end{eqnarray}
and \eqref{5.5} holds for some $\bar u\in C$ and $\bar v\in D$.
It suffices to show \eqref{5.4}.
We first show $l(u) \leq 0 $ for any $ u\in C$, i.e., $l\in -C^*$. Suppose to the contrary that there is $u^0\in C$ such that $l(u^0) > 0$. Since $X$ is  divisible and $C$ is a convex cone, from Proposition \ref{Pro:3.7} we have that $C=2^mC$ for any $m\in \mathbb{N}$. Then $2^mu^0\in C$ for any $m\in \mathbb{N}$. Since  $l:X\to \mathbb{R}$ is  additive, it follows that
\begin{eqnarray*}
l(2^mu^0)=2^ml(u^0)\to +\infty \quad {\rm as}\quad m\to+\infty,
 \end{eqnarray*}
a contradiction with \eqref{5.6}. Similarly, we can prove that $0 \leq l(v)$ for any $v\in D$, i.e., $l\in D^*$. Therefore, \eqref{5.4} holds. Since \eqref{5.4} and \eqref{5.5} are true, we have either$l(\bar u)<0$ or $l(\bar v)>0$. This completes the proof.
\cvd

%By Theorem \ref{Thm:5.4}, we have the following result by similar proof to that in Theorem \ref{Thm:5.6}.
%\begin{corollary}   \label{Cor:5.1} {\rm (Proper separation between a convex set and a convex cone)}  {\it Let $C$ be a nonempty convex subset and $D$ be a nonempty convex cone of a uniquely divisible commutative group $X$ with  $C\cap D = \emptyset$, and ${\rm prcore}\,(C-D)\neq\emptyset$  and $0\in D$. Then there is a nonzero additive function  $l:X\to \mathbb{R}$ such that
%\begin{eqnarray*}
%l(u) \leq 0 \leq l(v), \forall u\in C,\, \forall v\in D,
% \end{eqnarray*}
%or equivalently, $l\in (-C^*)\cap D^*$
%and
%$
%l(\bar u) <  l(\bar v)
%$
%for some $\bar u\in C$ and $\bar v\in D$, where either $l(\bar u)<0$ or $l(\bar v)>0$.
%}
%\end{corollary}

%{\bf Next theorem is a direct consequence of Theorem \ref{Thm:5.6}}

\begin{theorem} \label{Thm:5.7}  {\rm (Proper separation between  convex cones II)}  {\it Let $C$ and $D$ be nonempty subsets of a divisible commutative group $X$ with ${\rm cone}\,C\cap {\rm cone}\,D = \emptyset$ and ${\rm prcore}\,({\rm cone}\,C-{\rm cone}\,D)\neq\emptyset$.
Then there is a nonzero additive function  $l:X\to \mathbb{R}$  such that
\begin{eqnarray}\label{5.7}
l(u) \leq 0 \leq l(v), \forall u\in C,\, \forall v\in D,
 \end{eqnarray}
or equivalently, $l\in (-C^*)\cap D^*$ and
$
l(u^0) <  l(v^0)
$
for some $u^0\in C$ and $v^0\in D$, where either $l(u^0)<0$ or $l(v^0)>0$.
}
\end{theorem}
{\bf Proof.}\, From Theorem \ref{Thm:5.6} there is a nonzero additive function  $l:X\to \mathbb{R}$ such that
\begin{eqnarray*}
l(u) \leq 0\leq l(v), \forall u\in {\rm cone}\,C,\, \forall v\in {\rm cone}\,D,
 \end{eqnarray*}
or equivalently, $l\in (-({\rm cone}\,C)^*)\cap ({\rm cone}\,D)^*=(-C^*)\cap D^*$, in view of Proposition \ref{Pro:additiveope}, and
$
l(\bar v) <  l(\bar v)
$
holds for some $\bar u\in {\rm cone}\,C$ and $\bar v\in {\rm cone}\,D$ with either $l(\bar u)<0$ or $l(\bar v)>0$. As a consequence, \eqref{5.7} is true.
If $l(\bar u)<0$, then from Lemma \ref{LM:6.1} (ii) we have $l(u^0)<0$ some $u^0\in C$ and so $l(u^0)<l(v)$ for any $v\in D$. If  $l(\bar v)>0$, i.e., $-l(\bar v)<0$, then again from Lemma \ref{LM:6.1} (ii) we have $-l(v^0)<0$ some $v^0\in D$ and so $l(u)<l(v^0)$ for any $u\in C$. This proves that \eqref{5.7} holds and $l(u^0)<l(v^0)$ for some $u^0\in C$ and $v^0\in D$, where either $l(u^0)<0$ or $l(v^0)>0$.
\cvd

%Since $\bar u\in {\rm cone}\,C$ and $\bar v\in {\rm cone}\,D$,
%there are $m,n, r^i, s^j\in \mathbb{N}, x^i\in C, y^j\in D (i=1,\cdots, m, j=1,\cdots,n)$ such that
%$r\bar u=\sum_{i=1}^nr^i x^i$ and $s\bar v=\sum_{j=1}^ns^j y^j$. If $l(\bar u)<0$, then
%\begin{eqnarray*}
%\sum_{i=1}^n r^i l(x^i)=\sum_{i=1}^nl(r^i x^i)=l(\sum_{i=1}^nr^i x^i)=l(r\bar u)=rl(\bar u)<0,
%\end{eqnarray*}
%which implies that $l(x^{i_0})<0$ for some $x^{i_0}\in C$. If $l(\bar v)>0$, then by similar argument we can show $l(y^{j_0})>0$ for some $y^{j_0}\in D$.
%This, along with \eqref{5.7}, yields that $l(x^{i_0})< l(y^{j_0})$ holds for some $x^{i_0}\in C$ and $y^{j_0}\in D$  with either $l(x^{i_0})<0$  or $l(y^{j_0})>0$.
%This completes the proof.
%\cvd

\begin{theorem} \label{Thm:5.8}  {\rm (Proper separation between  convex cones III)}  {\it Let $C$ and $D$ be nonempty convex cones of a  divisible commutative group $X$ and ${\rm prcore}\,(C\backslash \{0\}-D)\neq\emptyset$, where
$C$ is pointed, i.e., $C\cap (-C)=\{0\}$, with $C\ne\{0\}$, such that $C\cap D = \{0\}$.
Then there is a nonzero additive function  $l:X\to \mathbb{R}$  such that \eqref{5.4} and \eqref{5.5} hold.
%\begin{eqnarray}\label{5.4}
%(u) \leq 0 \leq l(v), \forall u\in C,\, \forall v\in D,
% \end{eqnarray}
%or equivalently, $l\in (-C^*)\cap D^*$,
%and
%\begin{eqnarray}\label{5.5}
%l(\bar u) <  l(\bar v)
 %\end{eqnarray}
%for some $\bar u\in C$ and $\bar v\in D$, where either $l(\bar u)<0$ or $l(\bar v)>0$.
}
\end{theorem}
{\bf Proof.}\,  We first show that  $C\setminus\{0\}$ is a convex cone.  To, this aim, let $n, r^i\in \mathbb{N}, x, x^i\in C\backslash \{0\}  (i=1,\cdots, n,)$ such that
\begin{eqnarray}\label{5.8}
r x=\sum_{i=1}^nr^i x^i.
\end{eqnarray}
Since $C$ is a convex cone, one has $x\in C$. We declare that $x\not= 0$, i.e., $x\in C\setminus \{0\}$. Suppose to the contrary that $x=0$ in \eqref{5.8}.
If $n=1$, then from \eqref{5.8} we have $r^1x^1=0$, with $x^1\in C\setminus\{0\},$ which implies $r^1 \ge 2$. Then from Proposition \ref{Pro:3.6} (iii) we have $y:= (r^1-1)x^1\in C$ since $C$ is a convex cone, so that   $y+x^1=r^1x^1=0$ and
$y=-x^1\in -C\setminus\{0\}$, which contradicts the fact that $C$ is pointed.
Then without loss of generality,  we can suppose that in \eqref{5.8}, $n\ge 2$, with $r^ix^i \ne 0$, $i=1,\cdots,n$.  Again from Proposition \ref{Pro:3.6} (iii)  one has $z:= \sum_{i=1}^{n-1}r^i x^i\in C$ since $C$ is a convex cone, so that $z+r^n x^n=\sum_{i=1}^nr^ix^i=0$ and
$z=-r^nx^n\in -C\setminus\{0\}$, which contradicts also the fact that $C$ is pointed. As a consequence, $x\in C\setminus \{0\}$, which yields that $C\setminus \{0\}$ is a convex cone.
Since $C\cap D=\{0\}$ and $C$ is pointed,  $C\setminus \{0\}\cap D=\emptyset$.
Since $C\setminus\{0\}$ and $D$ are disjoint convex cones, by applying Theorem  \ref{Thm:5.6}, we obtain that  there exists a nonzero additive function  $l:X\to \mathbb{R}$  such that \eqref{5.4} and \eqref{5.5} hold,
recalling that $l(0)=0$.  \cvd

\begin{theorem} \label{Thm:5.9}  {\rm (Proper separation between $0$ and a convex cone IV)}  {\it Let $C$  be a nonempty convex cone of a   commutative group $X$ with   ${\rm qr core}\, C\not=\emptyset$ such that either $0\not\in C$ or
$C$ is not a subgroup of $X$ or
$C$ is pointed, i.e., $C\cap (-C)=\{0\}$, with $C\ne\{0\}$.
Then there is a nonzero additive function  $l:X\to \mathbb{R}$ such that $C$ and $0$ are properly separable, i.e.,
\begin{eqnarray}\label{5.9}
l(u) \leq l(0)=0, \forall u\in C,
 \end{eqnarray}
or equivalently, $l\in -C^*$, and
\begin{eqnarray}\label{5.10}
l(\bar u) <  0
 \end{eqnarray}
for some $\bar u\in C$.
}
\end{theorem}
{\bf Proof.}\, We preliminarly observe that $C={\rm cone}\,C$, by Proposition \ref{Pro:3.6} (i). If $0\not\in C$, then the thesis follows immediately from Theorem \ref{Thm:5.2}. If $C$ is not a subgroup of $X$, then
$0\not\in {\rm prcore}\,C$. If not, then $0 \in {\rm prcore}\,C$ and so $C={\rm cone}\,(C-0)$ is a subgroup of $X$, a contradiction. Consequently, the thesis follows immediately from Theorem \ref{Thm:5.3+}.
Assume that  $C$ is pointed with $C\ne\{0\}$. We observe that, setting $0\ne\bar x\in C$ allows that $-\bar x\not\in C$, since $C$ is pointed. This implies that $C$ is not a subgroup of $X$ and the thesis follows by the previous part of the proof.
\cvd
%As in the proof of the previous theorem, we can show that  $C^0:=C\setminus\{0\}$ is a convex cone.
%Since $C^0$ and $x^0:= 0 $ are disjoint,
%by applying Theorem  \ref{Thm:5.3}, we obtain that  there exists a nonzero additive function $l:X\to \mathbb{R}$  such that \eqref{5.9} and \eqref{5.10} hold,
%recalling that $l(0)=0$.

%{\bf Proof.}\, We preliminarly prove that $C-D$ is not a subgroup of $X$. By contradiction assume that $C-D$ is  a subgroup of $X$. In particuar let $0\ne c\in C$, we have that $c-0\in C-D$ implies that $-c\in C-D$. Then
%there exist $c_1\in C$, $d_1\in D$ such that $-c=c_1-d_1$. This implies that $d_1=c_1+c \in C$, since $C$ is a convex cone. By the assumption $C\cap D = \{0\}$ it must be $c_1+c= 0$, i.e., $c_1=-c$ which contradicts that
%$C$ is pointed.
%\cvd

\begin{remark} \label{Rm:5.2}{\it By applying Lemma \ref{LM:5.1} on the extension of additive functions, Theorems {\rm \ref{Thm:5.2}-\ref{Thm:5.4-}} prove new proper separation theorems between convex sets in commutative groups by removing the assumption that at least one of them has nonempty quasi core in {\rm \cite[Theorem 5.3]{LM24}}, where the Hahn-Banach theorem for groups {\rm \cite[Theorem 7]{BG18}} was employed. 
}
\end{remark}

\section{Applications to optimization problems with cone constraints in groups and semigroups}
\setcounter{equation}{0}
In this section, we investigate a class of optimization problems with cone constraints in groups and semigroups by exploiting the image space analysis \cite{Gia05}.
Without other specifications, we always assume that $X$ is a monoid and $\mathcal{D}$ is a  pointed  convex cone of a commutative group $Y$. %with $0\in \mathcal{D}$.

%{\bf Maybe we should use another letter different from $D$, to avoid misunderstanding with the set $D$ that appear in the Theorems of the previous section.}

Consider the following optimization problem with cone constraints (for short, OP):
\begin{eqnarray*}
 {\rm min}\,f(x)\quad {\rm subject\,\, to}\quad x\in X,\, g(x)\geq 0,
 \end{eqnarray*}
where $f: X \to \mathbb{R}$, $g: X\to Y$ and $\geq$ is a partial order defined by the  pointed convex cone $\mathcal{D}$, that is,
$x\geq y \Longleftrightarrow x-y\in \mathcal{D}, x, y\in Y$. From Proposition \ref{Pro:CPO}, one has $\geq$ is a compatible partial order.
Denote by $R^0:=\{x\in X: g(x)\geq 0\}$ the feasible set. We always assume that $R^0\neq \emptyset$.
% and $S$ the solution set of OP.

Observe that, $\bar x\in R^0$ is a solution of OP, if and only if
the system (in the unknown $x$):
\begin{eqnarray} \label{6.1}
f(\bar x)-f(x)>0,\, g(x)\geq 0,\, x\in X,
\end{eqnarray}
is impossible. Define the mapping $A_{\bar x}: X\to \mathbb{R}\times Y$ by
$
A_{\bar x}(x):=(f(\bar x)-f(x), g(x)).
$
Let
%We can associate CEP with the following sets:
\begin{eqnarray*}
{\mathcal{H}}&:=&\{(u,v)\in \mathbb{R}\times Y: u> 0,\, v\geq 0\}=\mathbb{R}_{++}\times \mathcal{D},\\
{\mathcal{K}}_{\bar x}&:=&\{(u,v)\in \mathbb{R}\times Y: u=f(\bar x)-f(x),\, v=g(x),\, x\in X\}=A_{\bar x}(X).
\end{eqnarray*}
 ${\mathcal{K}}_{\bar x}$ is called the image
associated with  OP, while  $\mathbb{R}\times Y$ is  the image space  associated with  OP.
%the system (\ref{5.1}).

Clearly, system (\ref{6.1}) is impossible, if and only if
\begin{eqnarray}\label{6.3}
{\mathcal{H}} \cap {\mathcal{K}}_{\bar x}=\emptyset.
\end{eqnarray}
Thus, $\bar x\in R^0$ is a solution of OP, if and only if (\ref{6.3}) holds.

In general, the image ${\mathcal{K}}_{\bar x}$  is not convex even when the
involved functions enjoy some convexity properties. To overcome this
difficulty,  we introduce a regularization of the image set  with respect to the cone $\mathbb{R}_+\times \mathcal{D}$, denoted by
$\mathcal{E}_{\bar x}:= {\mathcal{K}}_{\bar x}-(\mathbb{R}_+\times \mathcal{D})= A_{\bar x}(X)-(\mathbb{R}_+\times \mathcal{D})=\{(u,v)\in \mathbb{R}\times Y: u\leq f(\bar x)-f(x),\, v\le g(x),\, x\in X\}$,
which is called the extended image associated with OP. Since $(0,0)\in \mathbb{R}_+\times \mathcal{D}$, we have $\mathcal{E}_{\bar x}\supseteq{\mathcal{K}}_{\bar x}$. Since $R^0\neq \emptyset$, we have $(0,0)\in \mathcal{E}_{\bar x}$.
The following statement motivates the introduction of
the extended image.

\begin{proposition} \label{Pro:6.1}  {\rm  \cite[Proposition 5.1]{LM23}} {\it  The system {\rm(\ref{6.1})} is impossible, or {\rm(\ref{6.3})} holds, if and only if
\begin{eqnarray}\label{6.4}
{\mathcal{H}}\cap \mathcal{E}_{\bar x}=\emptyset
\end{eqnarray}
or, equivalently,
\begin{eqnarray}\label{6.4b}
{\mathcal{H}}_u\cap \mathcal{E}_{\bar x}=\emptyset,
\end{eqnarray}
where ${\mathcal{H}}_u:=\{(u,v)\in \mathbb{R}\times Y: u> 0,\, v = 0\}=\mathbb{R}_{++}\times \{0\}\subseteq {\mathcal{H}}$.
}
\end{proposition}

From Proposition \ref{Pro:6.1} we have the following.
\begin{proposition} \label{Pro:6.2} {\it  $\bar x\in R^0$ is a solution of OP, if and only if {\rm(\ref{6.4})}  or \eqref{6.4b} holds.
}
\end{proposition}

It is easy to check that $\mathcal{E}_{\bar x} -\mathcal{H}_u=\{(u,v)\in \mathbb{R}\times Y: u<f(\bar x)-f(x),\, v\le g(x),\, x\in X\}\subseteq \mathcal{E}_{\bar x}$. We have the following.

\begin{proposition} \label{Pro:6.3a}    {\it  Let $\bar x\in X$. Assume that $X$ is divisible and that $(Y, \leq)$ satisfies that for every $m\in \mathbb{N}, y^1, y^2\in Y$, $ my^1\leq m y^2\Longrightarrow y^1\leq y^2$. If  $f$ is convex and $g$ is  concave, then the set $\mathcal{E}_{\bar x} -\mathcal{H}_u$ is convex.
}
\end{proposition}
{\bf Proof.}\,  Let $n, r, r^i\in \mathbb{N}, (u^i, v^i)\in \mathcal{E}_{\bar x} -\mathcal{H}_u  (i=1,\cdots, n)$  such that
$r(u,v)=\sum_{i=1}^nr^i (u^i, v^i)$ with $\sum_{i=1}^n r^i=r$.
Then, there are $x^i\in X (i=1,\cdots, n)$ such that
\begin{eqnarray*}
 u^i < f(\bar x)- f(x^i),\quad v^i\le g(x^i), \quad i=1,\cdots, n.
\end{eqnarray*}
Summing the previous inequalities multiplied by $r^i$ yields that
\begin{eqnarray*}
\sum_{i=1}^n r^i u^i < \sum_{i=1}^n  r^i(f(\bar x)- f(x^i)),\quad \sum_{i=1}^n r^i v^i\le \sum_{i=1}^n r^i g(x^i).
\end{eqnarray*}
Since $X$ is divisible, there is $\hat x\in X$ such that $\sum_{i=1}^nr^i x^i=r\hat x$ and so, by the convexity of $f$ and the concavity of $g$ it follows that
\begin{eqnarray*}
r u= \sum_{i=1}^n r^i u^i < \sum_{i=1}^n  r^i(f(\bar x)- f(x^i))\le r(f(\bar x)-f(\hat x)),\quad rv=\sum_{i=1}^n r^i v^i\le \sum_{i=1}^n r^i g(x^i)\le rg(\hat x),
\end{eqnarray*}
which, by assumptions, implies that
$$ u < f(\bar x)- f(\hat x),\quad v\le g(\hat x), \quad i=1,\cdots, n,$$
i.e., $(u,v)\in \mathcal{E}_{\bar x} -\mathcal{H}_u$. \cvd

Under the convexity of $f$ and the concavity of $g$, we obtain the following  necessary optimality conditions for OP.
\begin{theorem} \label{Thm:6.1bb}  {\it  Let $\bar x\in R^0$ be an optimal solution of OP. Assume  that $X$ is divisible, $(Y, \leq)$ satisfies that for every $m\in \mathbb{N}, y^1, y^2\in Y$, $ my^1\leq m y^2\Longrightarrow y^1\leq y^2$ and that $f$ is convex and $g$ is concave. Assume that ${\rm prcore}\,(\mathcal{E}_{\bar x} -\mathcal{H}_u)\neq\emptyset$. Then   $(0,0)$ and $\mathcal{E}_{\bar x}$ are properly separable, that is, there is a nonzero additive function $l:=(\lambda, a)\in (\mathbb{R}_+)^*\times \mathcal{D}^*$  such that
\begin{eqnarray} \label{6.6-}
l(z)=\lambda(u)+a(v)\leq  0,\quad \forall z:=(u,v)\in \mathcal{E}_{\bar x}
\end{eqnarray}
and
\begin{eqnarray} \label{6.6-+}
l(\bar z)=\lambda(\bar u)+a(\bar v)< 0
\end{eqnarray}
for some $\bar z:=(\bar u,\bar v)\in \mathcal{E}_{\bar x}$.
}
\end{theorem}
{\bf Proof.}\,  By Proposition \ref{Pro:6.3a},   the set $\mathcal{E}_{\bar x} -\mathcal{H}_u$  is convex and since $\bar x$ is optimal for OP, by Proposition \ref{Pro:6.1}, we have that $(0,0)\not\in \mathcal{E}_{\bar x} -\mathcal{H}_u$. By applying Theorem \ref{Thm:5.2}  we have that  $\mathcal{E}_{\bar x} -\mathcal{H}_u$ and  $(0,0)$
 are properly linearly separable, i.e., there is a nonzero  additive function $l:=(\lambda, a)\in \mathbb{R}^*\times Y^*$ such that
\begin{eqnarray}\label{6.9+}
l(z) \leq l(0,0)=0,\, \forall z\in \,\mathcal{E}_{\bar x} -\mathcal{H}_u
 \end{eqnarray}
  and
\begin{eqnarray}\label{6.9b+}
l(\bar z)<0
 \end{eqnarray}
for some $\bar z\in \,\mathcal{E}_{\bar x} -\mathcal{H}_u \subseteq \mathcal{E}_{\bar x}$. Since $\lambda$ is additive, we have $\lambda(1)=\lambda(n\frac{1}{n})=n\lambda(\frac{1}{n})$ and so $\lambda(\frac{1}{n})=\frac{\lambda(1)}{n}$ for any $n\in \mathbb{N}$.
Then, from \eqref{6.9+} it follows that
\begin{eqnarray*}
\lambda(u-\frac{1}{n}) + a(v) \le 0,\quad \forall z=(u,v)\in \mathcal{E}_{\bar x},\quad \forall n\in \mathbb{N}
\end{eqnarray*}
and so
\begin{eqnarray*}
l(z)=\lambda(u) + a(v) \le \lambda(\frac{1}{n})=\frac{\lambda(1)}{n},\quad \forall z=(u,v)\in \mathcal{E}_{\bar x},\quad \forall n\in \mathbb{N}.
\end{eqnarray*}
Taking the limit in the previous inequality for $n \to +\infty$, yields that
\begin{eqnarray}\label{6.13N}
l(z)=\lambda(u) + a(v) \le 0,\quad \forall z=(u,v)\in \mathcal{E}_{\bar x},
\end{eqnarray}
which together with \eqref{6.9b+} proves that $(0,0)$ and $\mathcal{E}_{\bar x}$ are properly separable.

We now show that $l=(\lambda, a)\in (\mathbb{R}_+)^*\times \mathcal{D}^*$. Since $\mathcal{D}$ is a convex cone with $0\in \mathcal{D}$, from Proposition  \ref{Pro:3.6} (iii) one has $\mathcal{D}+\mathcal{D}=\mathcal{D}$.
Since $\bar x\in R^0$, $(0,0)\in  \mathcal{E}_{\bar x}$ and it follows that
\begin{eqnarray*}
-(0,d)=(0,0)-(0,d)\in A_{\bar x}(X)-(\mathbb{R}_+\times \mathcal{D})-(0\times \mathcal{D})=A_{\bar x}(X)-(\mathbb{R}_+\times \mathcal{D})=\mathcal{E}_{\bar x},\quad \forall d\in\mathcal{D}.
\end{eqnarray*}
Then from \eqref{6.13N} it follows that
\begin{eqnarray*}
-a(d)= a(-d) =\lambda (0)+a(-d)\le0, \forall d\in \mathcal{D},
\end{eqnarray*}
which proves that $a\in \mathcal{D}^*$.
Similarly,  from Proposition  \ref{Pro:3.6} (iii) and (iv) we have
\begin{eqnarray*}
-(u,0)=(0,0)-(u,0)\in A_{\bar x}(X)-(\mathbb{R}_+\times \mathcal{D})-(\mathbb{R}_+\times \mathcal{D})=A_{\bar x}(X)-(\mathbb{R}_+\times \mathcal{D})=\mathcal{E}_{\bar x},\quad \forall u\in\mathbb{R}_+.\end{eqnarray*}
Again from \eqref{6.13N} it follows that
\begin{eqnarray*}
-\lambda(u)= \lambda(-u)+ a(0) \le0, \quad \forall u\in \mathbb{R_+},
\end{eqnarray*}
which proves that $\lambda\in (\mathbb{R}_+)^*$.
This completes the proof.
\cvd

Observe that $\mathcal{E}_{\bar x}=A_{\bar x}(X)-(\mathbb{R}_+\times \mathcal{D})=(f(\bar x), 0)-[(f,-g)(X)+\mathbb{R}_+\times \mathcal{D}]$. The convexity of the extended image  $\mathcal{E}_{\bar x}$ can be guaranteed under convex-like assumptions.

\begin{proposition} \label{Pro:6.3}  {\rm  \cite[Proposition 5.3]{LM23}}  {\it  Let $\bar x\in X$. Then the following statements are true:
\begin{itemize}
\item[\rm(i)] Assume that $Y$ is divisible. If $\mathcal{E}_{\bar x}$ is convex, then $(f,-g)$ is convex-like.

\item[\rm(ii)] Assume that $(Y, \leq)$ satisfies that for every $m\in \mathbb{N}, y^1, y^2\in Y$, $ my^1\leq m y^2\Longrightarrow y^1\leq y^2$. If $(f,-g)$ is convex-like, then $\mathcal{E}_{\bar x}$ is convex.
\end{itemize}
}
\end{proposition}

\begin{proposition} \label{Pro:6.4a}  {\it  Let $\bar x\in X$. Assume that  $(Y, \leq)$ satisfies that for every $m\in \mathbb{N}, y^1, y^2\in Y$, $ my^1\leq m y^2\Longrightarrow y^1\leq y^2$ and $\mathcal{E}_{\bar x}$ is convex. Then
{\rm(\ref{6.4})} or equivalently \eqref{6.4b} holds if and only if
\begin{eqnarray}\label{6.5a}
{\mathcal{H}_u}\cap {\rm cone}\,\mathcal{E}_{\bar x}=\emptyset.
\end{eqnarray}
}
\end{proposition}
{\bf Proof.}\, If \eqref{6.5a} holds, then \eqref{6.4b} follows, since $\mathcal{E}_{\bar x}\subseteq {\rm cone}\,\mathcal{E}_{\bar x}$.

Assume that \eqref{6.4b} holds.
%From assumptions, Proposition \ref{Pro:6.3} yields that $\mathcal{E}_{\bar x}$ is convex.
%Since $\mathbb{R}_{++}$ and $\mathcal{D}$ are convex cones, from Proposition \ref{Pro:3.6} (iv) it follows that ${\mathcal{H}}$ is a convex cone.
Suppose to the contrary that {\rm(\ref{6.5a})} does not hold. Then there exists $(u,v)\in {\mathcal{H}_u}\cap {\rm cone}\,\mathcal{E}_{\bar x}$.
Then there are $n, r, r^i\in \mathbb{N}, (u^i, v^i)\in \mathcal{E}_{\bar x} (i=1,\cdots, n)$  such that
$r(u,v)=\sum_{i=1}^nr^i (u^i, v^i)$. In particular, we have
$$ ru= \sum_{i=1}^nr^i u^i,\quad rv= 0= \sum_{i=1}^nr^i  v^i,$$
noticing that $(u,v)\in \R_{++}\times\{0\}$.
Set $\bar u:=\frac{ru}{\sum_{i=1}^n {r^i}} \in\R_{++}$. Then $(\bar u,0)\in {\mathcal{H}_u}$ and
$$ (\sum_{i=1}^nr^i )(\bar u,0)= r(u,0)= \sum_{i=1}^nr^i (u^i, v^i)$$
which yields that $(\bar u,0)\in \mathcal{E}_{\bar x}$,
%Since $\mathbb{R}\times Y$ is divisible, there is $(u', v')\in \mathbb{R}\times Y$ such that $(\sum_{i=1}^nr^i)(u', v')=\sum_{i=1}^nr^i (u^i, v^i)$ and so $(u', v')\in \mathcal{E}_{\bar x}$,
since $\mathcal{E}_{\bar x}$ is convex.
% Consequently, $r(u,v)=(\sum_{i=1}^nr^i)(u', v')$.  Since $(u,v)\in {\mathcal{H}}$, from Proposition \ref{Pro:3.10} it follows that $r(u,v)=(\sum_{i=1}^nr^i)(u', v')\in {\mathcal{H}}$.
%Again by  Proposition \ref{Pro:3.10},  $(u', v')\in {\mathcal{H}}$
Thus $(\bar u, 0)\in {\mathcal{H}_u}\cap \mathcal{E}_{\bar x}$, a contradiction with \eqref{6.4b}.
\cvd

\begin{theorem} \label{Thm:6.1b}  {\it  Let $\bar x\in R^0$. Assume that ${\mathcal{H}}_u\cap{\rm cone}\,\mathcal{E}_{\bar x}=\emptyset$ and ${\rm prcore}\,({\rm cone}\,\mathcal{E}_{\bar x})\neq\emptyset$.  Then $\bar x$ is a solution of OP and  $(0,0)$ and $\mathcal{E}_{\bar x}$ are properly separable, that is, there is a nonzero additive function $l:=(\lambda, a)\in (\mathbb{R}_+)^*\times \mathcal{D}^*$  such that \eqref{6.6-} and \eqref{6.6-+} hold for some $\bar z:=(\bar u,\bar v)\in \mathcal{E}_{\bar x}$.
%\begin{eqnarray} \label{6.6-}
%\lambda(u)+a(v)=l(x)\leq  0,\quad \forall x:=(u,v)\in \mathcal{E}_{\bar x}
%\end{eqnarray}
%and
%\begin{eqnarray} \label{6.6-+}
%\lambda(\bar u)+a(\bar v)=l(\bar x)< 0
%\end{eqnarray}
%and $\inf_{(u,v)\in  \, {\mathcal{E}}_{\bar x}}(\lambda (u)+a(v))<0$.
}
\end{theorem}
{\bf Proof.}\,  Since  $\mathcal{E}_{\bar x}\subseteq {\rm cone}\,\mathcal{E}_{\bar x}$ and  ${\mathcal{H}}_u\cap{\rm cone}\,\mathcal{E}_{\bar x}=\emptyset$, we have  ${\mathcal{H}}_u\cap \mathcal{E}_{\bar x}=\emptyset$ and so from Propositions \ref{Pro:6.1} and \ref{Pro:6.2}   it follows that  $\bar x$ is a solution of OP. Let $z^0\in \mathcal{H}_u$. Then $z^0\not\in {\rm cone}\,\mathcal{E}_{\bar x}$. By applying Theorem \ref{Thm:5.2}  with $C:={\rm cone}\,\mathcal{E}_{\bar x}$ and $x^0:=z^0$, one has that there is a nonzero  additive function $l:=(\lambda, a)\in \mathbb{R}^*\times Y^*$ such that
\begin{eqnarray}\label{6.9a}
l(z) \leq l(z^0),\, \forall z\in {\rm cone}\,\mathcal{E}_{\bar x}
 \end{eqnarray}
  and
\begin{eqnarray}\label{6.9b}
 l(\hat{z}) <  l( z^0)
 \end{eqnarray}
for some $\hat{z}\in {\rm cone}\,\mathcal{E}_{\bar x}$. Since $(0,0)\in \mathcal{E}_{\bar x}\subseteq {\rm cone}\,\mathcal{E}_{\bar x}$ and $l(0,0)=0$, from \eqref{6.9a} we have $l(z^0)\ge 0$.
Since ${\rm cone}\,\mathcal{E}_{\bar x}$ is a convex cone, from Proposition \eqref{Pro:3.6} (iii) we have
$rz\in {\rm cone}\,\mathcal{E}_{\bar x}$ for any $r\in\mathbb{N}$ and  $z\in {\rm cone}\,\mathcal{E}_{\bar x}$. Since $l$  is additive,
from \eqref{6.9a}  it follows that
\begin{eqnarray*}
rl(z)=l(rz)\le l(z^0),\, \forall r\in\mathbb{N}, \,  \forall z\in {\rm cone}\,\mathcal{E}_{\bar x}
 \end{eqnarray*}
and so
\begin{eqnarray*}
l(z)\le \frac{l(z^0)}{r}\to 0\quad {\rm as}\,\, r\to +\infty, \  \forall z\in  {\rm cone}\,\mathcal{E}_{\bar x},
\end{eqnarray*}
which implies
\begin{eqnarray}\label{6.9aa}
l(z)\le 0=l(0,0), \  \forall z\in  \mathcal{E}_{\bar x}.
\end{eqnarray}

If $l(z^0) =0$ in \eqref{6.9b}, then,
by Lemma \ref{LM:6.1} (ii), there exists  some  $z^1 \in \mathcal{E}_{\bar x}$ such that
\begin{eqnarray*}
l(z^1) <  l(z^0)=0=l(0,0),
 \end{eqnarray*}
which, along with \eqref{6.9aa}, yields that   $(0, 0)$ and $\mathcal{E}_{\bar x}$ are properly separable.
Now assume  that $l(z^0) >0$. Then,  setting  $z^2:= (0,0) - z^0\in \mathcal{E}_{\bar x} - \mathcal{H}_u={\mathcal{K}}_{\bar x}-(\mathbb{R}_+\times \mathcal{D}+\mathcal{H}_u)={\mathcal{K}}_{\bar x}-\mathcal{H}\subseteq \mathcal{E}_{\bar x}$ and, since $l$ is additive,
 we have
 \begin{eqnarray*}
l(z^2)=l(-z^0)=-l(z^0)  < 0=l(0,0),
\end{eqnarray*}
which, along with \eqref{6.9aa}, yields that   $(0, 0)$ and $\mathcal{E}_{\bar x}$ are properly separable.

Similar to the proof in Theorem \ref{Thm:6.1bb}, we can show that $l=(\lambda, a)\in (\mathbb{R}_+)^*\times \mathcal{D}^*$. This completes the proof.
\cvd

\begin{remark} \label{Rem:6.2-} {\it We observe that, if ${\rm cone}\,\mathcal{E}_{\bar x}$ is pointed, then ${\mathcal{H}}_u\cap{\rm cone}\,\mathcal{E}_{\bar x}=\emptyset$. Suppose to the contrary that ${\mathcal{H}}_u\cap{\rm cone}\,\mathcal{E}_{\bar x}\neq\emptyset$. Let $(\bar u,0)\in{\mathcal{H}}_u\cap{\rm cone}\,\mathcal{E}_{\bar x}$, where $\bar u>0$. Note that ${\mathcal{H}}_u\subseteq \mathcal{H}$.
Since $\mathbb{R}_+, \mathbb{R}_{++}$ and $\mathcal{D}$ are convex cones, so are $\mathbb{R}_+\times \mathcal{D}$ and $\mathcal{H}$ by Proposition \ref{Pro:3.6} (iv). Since $(0,0)\in \mathbb{R}_+\times \mathcal{D}$,
from Proposition \ref{Pro:3.6} (iii) one has $\mathcal{H}=\mathbb{R}_+\times \mathcal{D}+\mathcal{H}$ and so
$\mathcal{E}_{\bar x}- \mathcal{H}= {\mathcal{K}}_{\bar x}-(\mathbb{R}_+\times \mathcal{D})-\mathcal{H}={\mathcal{K}}_{\bar x}-\mathcal{H}\subseteq \mathcal{E}_{\bar x}$.
Since  $(0,0)\in \mathcal{E}_{\bar x}$, it follows that $-(\bar u,0)= (0,0) - (\bar u,0)\in \mathcal{E}_{\bar x}\subseteq {\rm cone}\,\mathcal{E}_{\bar x}$, which is a contradiction with the fact that ${\rm cone}\, \mathcal{E}_{\bar x}$ is pointed.
}
\end{remark}

Under the assumption that $(f,-g)$ is convex-like, we obtain the following necessary optimality conditions for OP.
\vskip0.2cm
%{\bf If the previous proposition is true, then in the next theorem we can avoid the assumption that $Y$ is divisible (in fact, we can apply   Theorem \ref{Thm:6.1b} and Proposition \ref{Pro:6.6} (iv) to prove it. I have not done it yet since
%we should change the order of the results in the paper and I would like to have your opinion first)}

\begin{theorem} \label{Thm:6.1}  {\it  Let $\bar x\in R^0$. Assume that $(Y, \leq)$ satisfies that for every $m\in \mathbb{N}, y^1, y^2\in Y$, $ my^1\leq m y^2\Longrightarrow y^1\leq y^2$ and $(f,-g)$ is convex-like.
 Assume that ${\rm prcore}\,({\rm cone}\,\mathcal{E}_{\bar x})\neq\emptyset$.
If $\bar x\in R^0$ is a solution of OP, then   $(0,0)$ and $\mathcal{E}_{\bar x}$ are properly separable, that is, there is a nonzero additive function $l:=(\lambda, a)\in (\mathbb{R}_+)^*\times \mathcal{D}^*$ such that
\eqref{6.6-} and \eqref{6.6-+} hold for some $\bar z:=(\bar u,\bar v)\in \mathcal{E}_{\bar x}$, which yields that ${\mathcal{H}}$ and ${\mathcal{E}}_{\bar x}$ are properly separable.

}
\end{theorem}
{\bf Proof.}\,   By Proposition \ref{Pro:6.3} (ii), it follows that $\mathcal{E}_{\bar x}$ is convex. Since $\bar x\in R^0$ is a solution of OP, from Proposition  \ref{Pro:6.4a} we have that  \eqref{6.5a} holds.
Applying Theorem \ref{Thm:6.1b} we have that $(0,0)$ and $\mathcal{E}_{\bar x}$ are properly separable, that is, there is a nonzero additive function $l:=(\lambda, a)\in (\mathbb{R}_+)^*\times \mathcal{D}^*$ such that
\eqref{6.6-} and \eqref{6.6-+} hold for some $\bar z:=(\bar u,\bar v)\in \mathcal{E}_{\bar x}$. Since  $l=(\lambda, a)\in (\mathbb{R}_+)^*\times \mathcal{D}^*$, it follows that
\begin{eqnarray*}
0\leq l(z)=\lambda(u)+a(v) ,\quad \forall z:=(u,v)\in \mathcal{H},
\end{eqnarray*}
which, along with \eqref{6.6-} and \eqref{6.6-+} for some $\bar z:=(\bar u,\bar v)\in \mathcal{E}_{\bar x}$, yields that ${\mathcal{H}}$ and $\mathcal{E}_{\bar x}$ are properly separable,
This completes the proof.
\cvd

\begin{remark} \label{Rm:6.1}{\it Compared with {\rm  \cite[Theorem 5.3]{LM23}}, Theorem {\rm\ref{Thm:6.1}}  states the proper separation between ${\mathcal{H}}$ and $\mathcal{E}_{\bar x}$  by removing the assumption that $\mathcal{D}$ has nonempty quasi core and $Y$ is divisible. Moreover, from the proof in Theorem \ref{Thm:6.1b} we have $\inf_{(u,v)\in {\mathcal{H}}}(\lambda (u)+a(v))=\inf_{(u,v)\in \mathbb{R}_+\times \mathcal{D}}(\lambda (u)+a(v))=0$.

In fact, it suffices to show $\inf_{u\in \mathbb{R}_{++}}\lambda (u)=0$ since $l:=(\lambda, a)\in (\mathbb{R}_+)^*\times \mathcal{D}^*$ and $(0,0)\in \mathbb{R}_+ \times \mathcal{D}$. Suppose to the contrary that
$\eta:=\inf_{u\in \mathbb{R}_{++}}\lambda (u)>0$. Let $u^0\in \mathbb{R}_{++}$. Since $\lambda$ is additive, $\lambda(u^0)=m\lambda(\frac{u^0}{m})$ and so $\lambda(\frac{u^0}{m})=\frac{\lambda(u^0)}{m}\to 0$ as $m\to +\infty$, which is a contradiction with $\lambda(\frac{u^0}{m})\geq \eta>0$ for any $m\in \mathbb{N}$, since $\frac{u^0}{m}\in \mathbb{R}_{++}$ for any $m\in \mathbb{N}$.
This fills a small gap in the proof of necessity of {\rm  \cite[Proposition 5.6 (i) and (iii)]{LM23}}, which claims that $\lambda^0:=\lambda-\eta\in (\mathbb{R}_+)^*\times \mathcal{D}^*$ without proving $\eta=0$, however, $\lambda^0$ is additive if and only if $\eta=0$. The sufficiency follows immediately from the fact that if $\eta=0$, then $\lambda^0=\lambda\in (\mathbb{R}_+)^*$ and so $\lambda^0$ is additive. It suffices to show the necessity. Assume that $\lambda^0$ is additive. Since $\lambda\in (\mathbb{R}_+)^*$, one has $\eta=\lambda(0)-\lambda^0(0)=0-0=0$.
}
\end{remark}

In {\rm  \cite[Propositions 5.7 and 5.8]{LM23}}, Li and Mastroeni didn't prove (iv) $\Longleftrightarrow$ (iii) in the following proposition.

\begin{proposition} \label{Pro:6.5}  {\it  Let $\bar x\in X$. Then the following statements are true:
\begin{itemize}
\item[\rm(i)] The sets ${\mathcal{H}}$ and ${\mathcal{K}}_{\bar x}$ are  separable  if and only if there is a nonzero additive function $l:=(\lambda, a)\in (\mathbb{R}_+)^*\times \mathcal{D}^*$  such that
\begin{eqnarray}\label{6.6}
\sup_{(u,v)\in {\mathcal{K}}_{\bar x}}(\lambda (u)+a(v)) \leq 0.
 \end{eqnarray}

 \item[\rm(ii)] The sets ${\mathcal{H}}$ and ${\mathcal{E}}_{\bar x}$ are  separable if and only if there is a nonzero additive function $l:=(\lambda, a)\in (\mathbb{R}_+)^*\times \mathcal{D}^*$  such that
\begin{eqnarray}\label{6.7}
\sup_{(u,v)\in {\mathcal{E}}_{\bar x}}(\lambda (u)+a(v)) \leq 0.
 \end{eqnarray}

\item[\rm(iii)] The sets ${\mathcal{H}}$ and ${\mathcal{K}}_{\bar x}$ are  properly separable if and only if there is a nonzero  additive function $l:=(\lambda, a)\in (\mathbb{R}_+)^*\times \mathcal{D}^*$  such that \eqref{6.6}
holds and
$
\inf_{(u,v)\in {\mathcal{K}}_{\bar x}}(\lambda (u)+a(v)) < \sup_{(u,v)\in {\mathcal{H}}}(\lambda (u)+a(v)).
$

 \item[\rm(iv)] The sets ${\mathcal{H}}$ and ${\mathcal{E}}_{\bar x}$ are properly separable if and only if there is a nonzero additive function $l:=(\lambda, a)\in (\mathbb{R}_+)^*\times \mathcal{D}^*$  such that \eqref{6.7}
 holds and
$
\inf_{(u,v)\in {\mathcal{E}}_{\bar x}}(\lambda (u)+a(v)) < \sup_{(u,v)\in {\mathcal{H}}}(\lambda (u)+a(v)).
$

\end{itemize}
Moreover, {\rm (i) $\Longleftrightarrow$ (ii) $\Longleftarrow$ (iv) $\Longleftrightarrow$ (iii)}, where we assume that each statement holds for the same additive function $l$.

}
\end{proposition}
{\bf Proof.}\, It suffices to show  show (i) and (iii). Similarly, we can show (ii) and (iv).

(i) Necessity. Assume that the sets ${\mathcal{H}}$ and ${\mathcal{K}}_{\bar x}$ are  separable. Then there is a nonzero additive function $l:=(\lambda, a)\in \mathbb{R}^*\times Y^*$ such that
\begin{eqnarray*}
\sup_{(u,v)\in {\mathcal{K}}_{\bar x}}(\lambda (u)+a(v)) \leq \inf_{(u,v)\in {\mathcal{H}}}(\lambda (u)+a(v)).
 \end{eqnarray*}
From Remark \ref{Rm:6.1} we have $\inf_{(u,v)\in {\mathcal{H}}}(\lambda (u)+a(v))=\inf_{(u,v)\in \mathbb{R}_+\times \mathcal{D}}(\lambda (u)+a(v))=0$ and $l:=(\lambda, a)\in (\mathbb{R}_+)^*\times \mathcal{D}^*$.
This proves \eqref{6.6} holds.

Sufficiency. Since $(\lambda, a)\in (\mathbb{R}_+)^* \times \mathcal{D}^*$, we have $\inf_{(u,v)\in {\mathcal{H}}}(\lambda(u)+a(v))\geq 0$ and from \eqref{6.6} it follows that
$
\sup_{(u,v)\in {\mathcal{K}}_{\bar x}}(\lambda (u)+a(v)) \leq \inf_{(u,v)\in {\mathcal{H}}}(\lambda (u)+a(v)).
$
This yields that ${\mathcal{H}}$ and ${\mathcal{K}}_{\bar x}$ are  separable.

(iii) Necessity. Assume that the sets ${\mathcal{H}}$ and ${\mathcal{K}}_{\bar x}$ are properly separable. Then there is a nonzero additive function $l:=(\lambda, a)\in \mathbb{R}^*\times Y^*$ such that
\eqref{6.6} holds and
$
\inf_{(u,v)\in {\mathcal{K}}_{\bar x}}(\lambda(u)+a(v)) < \sup_{(u,v)\in {\mathcal{H}}}(\lambda(u)+a(v)).
$
From the proof of the Necessity in (i), we have $(\lambda, a)\in  (\mathbb{R}_+)^* \times \mathcal{D}^*$.

Sufficiency. It is clear.

Observe that (ii) $\Longleftarrow$ (iv) $\Longleftarrow$ (iii) are clear.
We first show (i) $\Longleftrightarrow$ (ii). It suffices to show (i) $\Longrightarrow$ (ii) since (i) $\Longleftarrow$ (ii) is clear. To this aim, suppose that (i) holds, i.e., there is a nonzero additive function $l:=(\lambda, a)\in (\mathbb{R}_+)^*\times \mathcal{D}^*$  such that \eqref{6.6} holds. For any $(u,v)\in {\mathcal{E}}_{\bar x}$, there is $(u^0,v^0)\in {\mathcal{H}}$ such that $(u,v)+(u^0,v^0)\in {\mathcal{K}}_{\bar x}$.  From Remark \ref{Rm:6.1}  we have $\inf_{(u,v)\in {\mathcal{H}}}(\lambda (u)+a(v))=0$.
Since $(u,v)=((u,v)+(u^0,v^0))-(u^0,v^0)$ and $l:=(\lambda, a)\in (\mathbb{R}_+)^*\times \mathcal{D}^*$ is additive, it follows that
\begin{eqnarray}\label{6.8-}
(\lambda(u)+a(v))&=&(\lambda(u+u^0)+a(v+v^0))-(\lambda(u^0)+a(v^0))\nonumber\\
                 &\leq& 0,
\end{eqnarray}
where the inequality follows from \eqref{6.6} and the fact that $\inf_{(u,v)\in {\mathcal{H}}}(\lambda (u)+a(v))=0$. Since $(u,v)\in {\mathcal{E}}_{\bar x}$ is arbitrary, from \eqref{6.8-} it follows that \eqref{6.7} holds, proving (ii).

We now show (iv) $\Longrightarrow$ (iii).
Suppose that there is a nonempty additive function $l:=(\lambda, a)\in (\mathbb{R}_+)^*\times \mathcal{D}^*$  such that \eqref{6.7} holds and $\inf_{(u,v)\in {\mathcal{E}}_{\bar x}}(\lambda (u)+a(v)) < \sup_{(u,v)\in {\mathcal{H}}}(\lambda (u)+a(v))$. Since  $\sup_{(u,v)\in {\mathcal{H}}}(\lambda (u)+a(v))\geq \inf_{(u,v)\in {\mathcal{H}}}(\lambda (u)+a(v))=0$, we first declare that either $\sup_{(u,v)\in {\mathcal{H}}}(\lambda (u)+a(v))=0$ or $\sup_{(u,v)\in {\mathcal{H}}}(\lambda (u)+a(v))=+\infty$. In fact, if $\sup_{(u,v)\in {\mathcal{H}}}(\lambda (u)+a(v))\neq0$, i.e., $\sup_{(u,v)\in {\mathcal{H}}}(\lambda (u)+a(v))>0$, then there is $(u^0,v^0)\in {\mathcal{H}}$ such that
$\lambda (u^0)+a(v^0)>0$.
Note that ${\mathcal{H}}$ is a convex cone. From Proposition \ref{Pro:3.6} (iii) we have $m(u^0,v^0)\in {\mathcal{H}}$ for any $m\in \mathbb{N}$ and so $\lambda (mu^0)=m\lambda (u^0)$ and $a(mv^0)=ma(v^0)$.
Thus, $0<\lambda (mu^0)+a(mv^0)=m(\lambda (u^0)+a(v^0))\to +\infty$ as $m\to +\infty$, which yields that $\sup_{(u,v)\in {\mathcal{H}}}(\lambda (u)+a(v))=+\infty$.
Since $\inf_{(u,v)\in {\mathcal{E}}_{\bar x}}(\lambda (u)+a(v)) < \sup_{(u,v)\in {\mathcal{H}}}(\lambda (u)+a(v))$, there is $(\bar u, \bar v) \in {\mathcal{E}}_{\bar x}$ such that $\lambda (\bar u)+a(\bar v)< \sup_{(u,v)\in {\mathcal{H}}}(\lambda (u)+a(v))$. Then there are $(\mathbf{u},\mathbf{v})\in {\mathcal{K}}_{\bar x}$ and $(u^0, v^0)\in {\mathcal{H}}$ such that
$(\bar u, \bar v)=(\mathbf{u},  \mathbf{v})-(u^0, v^0)$ and so $(\mathbf{u},  \mathbf{v})=(\bar u, \bar v)+(u^0, v^0)$. As a consequence,
\begin{eqnarray*}
\inf_{(u,v)\in {\mathcal{K}}_{\bar x}}(\lambda (u)+a(v))&\leq&\lambda(\mathbf{u})+a(\mathbf{v})\\
&=&(\lambda(\bar u)+a(\bar v))+(\lambda(u^0)+a(v^0))\\
&<& \sup_{(u,v)\in {\mathcal{H}}}(\lambda (u)+a(v))+(\lambda(u^0)+a(v^0))\\
&\leq& \sup_{(u,v)\in {\mathcal{H}}}(\lambda (u)+a(v))+ \sup_{(u,v)\in {\mathcal{H}}}(\lambda (u)+a(v))\\
&=&\sup_{(u,v)\in {\mathcal{H}}}(\lambda (u)+a(v)),
\end{eqnarray*}
where the second strict inequality follows from the fact that $\lambda (\bar u)+a(\bar v)< \sup_{(u,v)\in {\mathcal{H}}}(\lambda (u)+a(v))$, and the last equality follows from the fact that either $\sup_{(u,v)\in {\mathcal{H}}}(\lambda (u)+a(v))=0$ or $\sup_{(u,v)\in {\mathcal{H}}}(\lambda (u)+a(v))=+\infty$.
This completes the proof.
\cvd

We have the following equivalent characterizations of  proper separation.

\begin{proposition} \label{Pro:6.6}  {\it  Let $\bar x\in X$. Then the following statements are true:
\begin{itemize}
\item[\rm(i)]  $(0,0)$ and ${\mathcal{E}}_{\bar x}$ are  separable if and only if there is a nonzero additive function $l:=(\lambda, a)\in (\mathbb{R}_+)^*\times \mathcal{D}^*$  such that \eqref{6.7} holds.

\item[\rm(ii)]  $(0,0)$ and ${\mathcal{E}}_{\bar x}$ are properly  separable if and only if there is a nonzero additive function $l:=(\lambda, a)\in (\mathbb{R}_+)^*\times \mathcal{D}^*$  such that \eqref{6.7} holds and $\inf_{(u,v)\in {\mathcal{E}}_{\bar x}}(\lambda (u)+a(v))<0$.

\item[\rm(iii)] The sets ${\mathcal{H}}$ and ${\mathcal{E}}_{\bar x}$ are  separable if and only if $(0,0)$ and ${\mathcal{E}}_{\bar x}$ are  separable.

\item[\rm(iv)] The sets ${\mathcal{H}}$ and ${\mathcal{E}}_{\bar x}$ are properly separable if and only if $(0,0)$ and ${\mathcal{E}}_{\bar x}$ are properly separable.
\end{itemize}
}
\end{proposition}
{\bf Proof.}\, (i) Sufficiency. It is clear.

Necessity. Assume that $(0,0)$ and ${\mathcal{E}}_{\bar x}$ are  separable. Then  there is a nonzero additive function $l:=(\lambda, a)\in \mathbb{R}^*\times Y^*$  such that \eqref{6.7} holds.
Since $\mathbb{R}_+$ and $\mathcal{D}$ are convex cones, so is $\mathbb{R}_+\times \mathcal{D}$ by Proposition \ref{Pro:3.6} (iv). Since $(0,0)\in \mathbb{R}_+\times \mathcal{D}$,
Proposition \ref{Pro:3.6} (iii) implies that $\mathbb{R}_+\times \mathcal{D}=\mathbb{R}_+\times \mathcal{D}+\mathbb{R}_+\times \mathcal{D}$. As a consequence,
$\mathcal{E}_{\bar x}= {\mathcal{K}}_{\bar x}-(\mathbb{R}_+\times \mathcal{D})={\mathcal{K}}_{\bar x}-(\mathbb{R}_+\times \mathcal{D}+\mathbb{R}_+\times \mathcal{D})=\mathcal{E}_{\bar x}-\mathbb{R}_+\times \mathcal{D}$.
Now \eqref{6.7} yields that
\begin{eqnarray*}
\sup_{(u,v)\in {\mathcal{E}}_{\bar x}}(\lambda (u)+a(v)) \leq \inf_{(u,v)\in \mathbb{R}_+\times \mathcal{D}}(\lambda (u)+a(v))=0
 \end{eqnarray*}
and so $l:=(\lambda, a)\in (\mathbb{R}_+)^*\times \mathcal{D}^*$, where the equality follows from Remark \ref{Rm:6.1}.

(ii) The proof is similar to that in (i).

(iii) It follows immediately from (i) and Proposition \ref{Pro:6.5} (ii).

(iv) Sufficiency. Assume that $(0,0)$ and ${\mathcal{E}}_{\bar x}$ are properly  separable. Then from (ii) we have there is a nonzero additive function $l:=(\lambda, a)\in (\mathbb{R}_+)^*\times \mathcal{D}^*$  such that \eqref{6.7} holds and $\inf_{(u,v)\in {\mathcal{E}}_{\bar x}}(\lambda (u)+a(v))<0$. Since $0\leq\sup_{(u,v)\in {\mathcal{H}}}(\lambda (u)+a(v))$, one has $\inf_{(u,v)\in {\mathcal{E}}_{\bar x}}(\lambda (u)+a(v))< \sup_{(u,v)\in {\mathcal{H}}}(\lambda (u)+a(v))$. By Proposition \ref{Pro:6.5} (iv) it follows that the sets ${\mathcal{H}}$ and ${\mathcal{E}}_{\bar x}$ are properly separable.

Necessity. Assume that ${\mathcal{H}}$ and $\mathcal{E}_{\bar x}$ are properly separable. Then from Proposition \ref{Pro:6.5} (iv), there exists a nonzero additive function $l:=(\lambda, a)\in (\mathbb{R}_+)^* \times \mathcal{D}^*$ such that \eqref{6.7} and $\inf_{\mathbf{x}\in {\mathcal{E}}_{\bar x}}l(\mathbf{x})=\inf_{(u,v)\in {\mathcal{E}}_{\bar x}}(\lambda (u)+a(v))< \sup_{(u,v)\in {\mathcal{H}}}(\lambda (u)+a(v))=\sup_{\mathbf{y}\in {\mathcal{H}}}l(\mathbf{y})$ hold. Then there are $\bar{\mathbf{x}}\in {\mathcal{E}}_{\bar x}$ and $\bar{\mathbf{y}}\in {\mathcal{H}}$ such that $l(\bar{\mathbf{x}})<l(\bar{\mathbf{y}})$
and so $l(\bar{\mathbf{x}}-\bar{\mathbf{y}})<0$, since $l$ is  additive. Then $\bar{\mathbf{z}}:=\bar{\mathbf{x}}-\bar{\mathbf{y}}\in \mathcal{E}_{\bar x} -\mathcal{H} \subseteq\mathcal{E}_{\bar x}$ (see Remark \ref{Rem:6.2-})
and consequently, $l(\bar{\mathbf{z}})<0$, which along with \eqref{6.7}, yields that $(0,0)$  and $\mathcal{E}_{\bar x}$ are properly separable.
\cvd

\begin{definition} \label{Deg:lanrange}  {\it The Lagrangian function $L:(\mathbb{R}_+)^* \times \mathcal{D}^*\times X\to \mathbb{R}$ associated with OP is
defined by $L(\lambda,a,x):= \lambda (f(x)) -a(g(x)), (\lambda,a,x)\in (\mathbb{R}_+)^* \times \mathcal{D}^*\times X$. Let $\lambda^*\in (\mathbb{R}_+)^* $.
The point $(a^*,\bar x)\in \mathcal{D}^*\times X$ is said to be a saddle point for $L(\lambda^*,\cdot,\cdot)$ on  $ \mathcal{D}^*\times X$,  if
 \begin{eqnarray*}
L(\lambda^*,a,\bar x)\le L(\lambda^*,a^*,\bar x)\le L(\lambda^*,a^*, x),  \quad\forall (a,x)\in  \mathcal{D}^*\times X.
\end{eqnarray*}
 }
\end{definition}

Consider the following Slater-type  constraint qualification  for OP:
\begin{eqnarray}\label{CQ}
0\in{\rm prcore (cone}\, (g(X)-\mathcal{D})).
\end{eqnarray}

In the following, we obtain saddle points for the Lagrangian function associated with OP under the   Slater-type constraint qualification \eqref{CQ}.

\begin{theorem} \label{Th:6.2}  {\it Let $\bar x\in R^0$. Assume that  ${\mathcal{H}}_u\cap{\rm cone}\,\mathcal{E}_{\bar x}=\emptyset$  with ${\rm prcore}\,({\rm cone}\,\mathcal{E}_{\bar x})\neq\emptyset$  and that Slater-type constraint qualification \eqref{CQ} holds. Then  $\bar x$ is a solution of OP and there exists a nonzero additive function $l^*:=(\lambda^*, a^*)\in (\mathbb{R}_+)^* \times \mathcal{D}^*$, with $\lambda^*\not\equiv 0$,  such that $(a^*,\bar x)$ is a saddle point for
$L(\lambda^*,\cdot,\cdot)$  on $\mathcal{D}^*\times X$.
 }
\end{theorem}
{\bf Proof.} From Theorem \ref{Thm:6.1b} we have $\bar x$ is a solution of OP and  $(0,0)$  and $\mathcal{E}_{\bar x}$ are properly separable. That is, there exists a nonzero additive function $l^*:=(\lambda^*, a^*)\in (\mathbb{R}_+)^* \times \mathcal{D}^*$ such that
\begin{equation}\label{6.8}
\lambda^*(u) + a^*(v)=l^*(u,v)\le l^*(0,0), \quad \forall (u,v)\in\mathcal{E}_{\bar x}
\end{equation}
and  the strict inequality in \eqref{6.8} holds for  some $(\hat u,\hat v)\in \mathcal{E}_{\bar x}$.
We now show that $\lambda^*\not\equiv 0$. Suppose to the contrary that $\lambda^*\equiv 0$. Then from \eqref{6.8} we have
\begin{equation}\label{6.9}
a^*(v) \le 0, \quad \forall  (u,v)\in \mathcal{E}_{\bar x}
\end{equation}
and  the strict inequality in \eqref{6.9} holds for  some $(\hat u,\hat v)\in  \mathcal{E}_{\bar x}$.
Observe that for any $v^0\in g(X)-\mathcal{D}$, there are $x^0\in X$ and $d^0\in \mathcal{D}$ such that $v^0=g(x^0)-d^0$. Setting $(u^0,v^0):=(f(\bar x)-f(x^0), g(x^0))-(0,d^0)$ allows that $(u^0,v^0)\in \mathcal{E}_{\bar x}$.
This, along with \eqref{6.9}, yields that
\begin{equation}\label{6.10}
a^*(v) \le 0, \quad \forall v\in g(X)-\mathcal{D}
\end{equation}
and  by Lemma \ref{LM:6.1} (ii) it follows that
\begin{equation}\label{6.11}
a^*(\hat v) < 0
\end{equation}
for some  $(\hat u,\hat v):=(f(\bar x)-f(\hat x), g(\hat x))-(\hat t, \hat d)\in  \mathcal{E}_{\bar x}$, where $\hat x\in X$ and $(\hat t, \hat d)\in \mathbb{R}_+\times \mathcal{D}$.
That is, \eqref{6.11} holds for some $\hat v\in g(X)-\mathcal{D}$.
Again, by Lemma \ref{LM:6.1} (i) and (ii), it follows from  \eqref{6.10} and \eqref{6.11} that
\begin{equation}\label{6.12}
a^*(v) \le 0, \quad \forall v\in {\rm cone}\,(g(X)-\mathcal{D})
\end{equation}
and
\begin{equation}\label{6.13}
a^*(\hat v) < 0, \quad {\rm for \, some}\,\, \hat v\in {\rm cone}\,(g(X)-\mathcal{D}).
\end{equation}
Since the Slater-type constraint qualification  \eqref{CQ} holds, $0\in {\rm cone}\,(g(X)-\mathcal{D})$ and ${\rm cone}({\rm cone}\,(g(X)-\mathcal{D})-0)={\rm cone}\,(g(X)-\mathcal{D})$ is a subgroup of $Y$, where the equality follows from Proposition \ref{Pro:3.6} (i). As a consequence, $\pm v\in {\rm cone}\,(g(X)-\mathcal{D})$  for any $v\in {\rm cone}\,(g(X)-D)$. Since $a^*$ is additive, from \eqref{6.12} it follows that
$a^*(v) = 0$ for any $v\in {\rm cone}\,(g(X)-\mathcal{D})$,
which contradicts  \eqref{6.13}.
This proves that  $\lambda^*\not\equiv 0$.

Since $\mathcal{K}_{\bar x}\subseteq  \mathcal{E}_{\bar x}$, \eqref{6.8}  implies that
\begin{equation}\label{6.14}
\lambda^*(f(\bar x)-f(x)) + a^*(g(x)) \le 0 \quad \forall x\in X.
\end{equation}
Setting $x:=\bar x$ in \eqref{6.14} yields that $a^*(g(\bar x))\leq 0$. Since $a^*\in \mathcal{D}^*$ and  $\bar x\in R^0$, i.e., $g(\bar x)\in \mathcal{D}$, one has $a^*(g(\bar x))\ge0$ and
consequently,
\begin{equation}\label{6.15}
a^*(g(\bar x))=0,
\end{equation}
which implies that
$-a(g(\bar x)) \le -a^*(g(\bar x))=0$ for any $a\in \mathcal{D}^*$. As a consequence,
\begin{equation*}
L(\lambda^*,a,\bar x)=\lambda^*(f(\bar x))-a(g(\bar x)) \le \lambda^*(f(\bar x)) -a^*(g(\bar x))= L(\lambda^*,a^*,\bar x),\quad\forall a\in \mathcal{D}^*.
\end{equation*}
From  \eqref{6.14} and \eqref{6.15} we also have
\begin{equation*}
L(\lambda^*,a^*,\bar x)=\lambda^*(f(\bar x)) -a^*(g(\bar x))\le \lambda^*(f(x)) -a^*(g(x))=L(\lambda^*,a^*, x),\quad\forall x\in X.
\end{equation*}
This proves that $(a^*,\bar x)$ is a saddle point for $L(\lambda^*,\cdot,\cdot)$ on $ \mathcal{D}^*\times X$.
\cvd

%{\color{red}

\begin{remark}\label{Rem:6.2} {\it
%We observe that if ${\rm cone}\, \mathcal{E}_{\bar x}$ is pointed, then $\bar x\in R^0$ is a solution of OP. In fact, suppose to the contrary  that  $\bar x$ is not a solution of OP. Then from Proposition \ref{Pro:6.2} one has %\eqref{6.4} does not hold and so there exists $(\bar u,\bar v)\in \mathcal{E}_{\bar x}\cap \mathcal{H}\subseteq {\rm cone}\,\mathcal{E}_{\bar x}\setminus\{(0,0)\}$. Since $\mathbb{R}_+, \mathbb{R}_{++}$ and $\mathcal{D}$ are %convex cones, so are $\mathbb{R}_+\times \mathcal{D}$ and $\mathcal{H}$ by Proposition \ref{Pro:3.6} (iv). Since $(0,0)\in \mathbb{R}_+\times \mathcal{D}$,
%from Proposition \ref{Pro:3.6} (iii) one has $\mathcal{H}=\mathbb{R}_+\times \mathcal{D}+\mathcal{H}$ and so
%$\mathcal{E}_{\bar x}- \mathcal{H}= {\mathcal{K}}_{\bar x}-(\mathbb{R}_+\times \mathcal{D})-\mathcal{H}={\mathcal{K}}_{\bar x}-\mathcal{H}\subseteq \mathcal{E}_{\bar x}$.
%Since  $(0,0)\in \mathcal{E}_{\bar x}$, it follows that $-(\bar u,\bar v)= (0,0) - (\bar u,\bar v)\in \mathcal{E}_{\bar x}\subseteq {\rm cone}\,\mathcal{E}_{\bar x}$, which is a contradiction with the fact that ${\rm cone}\, %\mathcal{E}_{\bar x}$ is pointed.
Similarly to the previous remark,  in Theorem \ref{Th:6.2}, we may simply assume that ${\rm cone}\,\mathcal{E}_{\bar x}$ is pointed, which implies ${\mathcal{H}}_u\cap{\rm cone}\,\mathcal{E}_{\bar x}=\emptyset$. In this case, since  $(0,0)\in \mathcal{E}_{\bar x}\subseteq {\rm cone}\,\mathcal{E}_{\bar x}$, we have
${\rm cone}\,({\rm cone}\,\mathcal{E}_{\bar x}-(0,0))={\rm cone}\,\mathcal{E}_{\bar x}$ is not a subgroup of $\mathbb{R}\times Y$, that is, $(0,0)\not\in{\rm prcore (cone}\,{\mathcal{E}}_{\bar x})$, which implies that
$(0,0)\not\in{\rm qcore (cone}\,{\mathcal{E}}_{\bar x})$. However, in this way we avoid to assume  that $Y$ is divisible, compared with {\rm \cite[Theorems 6.2 and 6.3]{LM25}}.
}
\end{remark}

In the following result, we remove the assumption of nonemptyness of the quasi core of  ${\rm cone}\, \mathcal{E}_{\bar x}$, compared with {\rm \cite[Theorem 6.3]{LM25}}.

\begin{theorem} \label{Th:6.3}  {\it Let $\bar x\in R^0$ be a solution of OP. Assume  that   $(Y, \leq)$ satisfies that for every $m\in \mathbb{N}, y^1, y^2\in Y$, $ my^1\leq m y^2\Longrightarrow y^1\leq y^2$, $(f,-g)$ is convex-like, ${\rm prcore}\,({\rm cone}\,\mathcal{E}_{\bar x})\neq\emptyset$ and the  Slater-type constraint qualification \eqref{CQ} holds.  Then  there exists a nonzero additive function $l^*:=(\lambda^*, a^*)\in (\mathbb{R}_+)^* \times \mathcal{D}^*$, with $\lambda^*\not\equiv 0$,  such that $(a^*,\bar x)$ is a saddle point for $L(\lambda^*,\cdot,\cdot)$  on $\mathcal{D}^*\times X$.
 }
\end{theorem}
{\bf Proof.} Since $\bar x\in R^0$ is a solution of OP, then from Propositions \ref{Pro:6.2}   and  \ref{Pro:6.4a}, it follows that ${\mathcal{H}}_u\cap{\rm cone}\,\mathcal{E}_{\bar x}=\emptyset$ and  Theorem \ref{Thm:6.1b} yields that $(0,0)$  and $\mathcal{E}_{\bar x}$ are properly separable. The subsequent part of the  proof is similar to that in Theorem \ref{Th:6.2}, so we omit it.
\cvd

%From  Theorem \ref{Thm:6.1} we have ${\mathcal{H}}$ and $\mathcal{E}_{\bar x}$ are properly separable, i.e., there exists a nonzero additive function $l^*:=(\lambda^*, a^*)\in (\mathbb{R}_+)^* \times \mathcal{D}^*$  such that  %\eqref{6.6-} and \eqref{6.6-+} hold for some $\bar x:=(\bar u,\bar v)\in \mathcal{E}_{\bar x}$ and $\bar y:=(\bar{\mathbf{u}},\bar{\mathbf{v}})\in {\mathcal{H}}$,  where either $l^*(\bar x)<0$ or $l^*(\bar y)>0$. Now, %Proposition \ref{Pro:6.6} (iv) yields that $(0,0)$  and $\mathcal{E}_{\bar x}$ are properly separable by $l^*$.
%That is,
%\begin{equation}\label{6.16}
%\lambda^*(u) + a^*(v)=l^*(u,v)\le l^*(0,0)=\lambda^*(0) + a^*(0)=0, \quad \forall (u,v)\in  \mathcal{E}_{\bar x}
%\end{equation}
%and the strict inequality in \eqref{6.16} holds for  some $(\hat u,\hat v)\in  \mathcal{E}_{\bar x}$.
%We now show that $\lambda^*\not\equiv 0$. Suppose to the contrary that $\lambda^*\equiv 0$. Then from \eqref{6.16} we have
%\begin{equation}\label{6.17}
%a^*(v) \le 0, \quad \forall  (u,v)\in  \mathcal{E}_{\bar x},
%\end{equation}
%and the strict inequality in \eqref{6.17} holds for  some $(\hat u,\hat v)\in  \mathcal{E}_{\bar x}$.
%Observe that for any $v^0\in g(X)-\mathcal{D}$, there are $x^0\in X$ and $d^0\in \mathcal{D}$ such that $v^0=g(x^0)-d^0$. Setting $(u^0,v^0):=(f(\bar x)-f(x^0), g(x^0))-(0,d^0)$ yields that $(u^0,v^0)\in \mathcal{E}_{\bar x}$.
%The subsequent part of the  proof is similar to that in Theorem \ref{Th:6.2}, so we omit it.
%\cvd

Next result can similarly be proved by Theorem  \ref{Thm:6.1bb} or, alternatively, it  follows from the previous theorem, taking into account Proposition \ref{Pro:4.5}.

\begin{theorem} \label{Thm:6.4}  {\it  Let $\bar x\in R^0$ be an optimal solution of OP. Assume  that $X$ is divisible, $(Y, \leq)$ satisfies that for every $m\in \mathbb{N}, y^1, y^2\in Y$, $ my^1\leq m y^2\Longrightarrow y^1\leq y^2$, $f$ is convex and $g$ is concave,  ${\rm prcore}\,({\rm cone}\,\mathcal{E}_{\bar x})\neq\emptyset$  and  the  Slater-type constraint qualification \eqref{CQ} holds. Then  there exists a nonzero additive function $l^*:=(\lambda^*, a^*)\in (\mathbb{R}_+)^* \times \mathcal{D}^*$, with $\lambda^*\not\equiv 0$,  such that $(a^*,\bar x)$ is a saddle point for $L(\lambda^*,\cdot,\cdot)$  on $\mathcal{D}^*\times X$.
}
\end{theorem}
{\bf Proof.}  Since $f$ is convex and $g$ is concave, $(f,-g)$ is convex. Since $X$ divisible, from Proposition  \ref{Pro:4.5} it follows that $(f,-g)$ is convex-like.
Now, applying Theorem \ref{Th:6.3} yields the desired conclusion.
\cvd

The following result relies on the properties of quasi-relative core of the set $\mathcal{E}_{\bar x} $ in order to achieve separation and the existence of a saddle point for $L$.

\begin{theorem} \label{Thm:6.7}  {\it  Let $\bar x\in R^0$. Assume  that  $(0,0)\not\in {\rm prcore}\,\mathcal{E}_{\bar x}\neq\emptyset$   and  the  Slater-type constraint qualification \eqref{CQ} holds. Then  there exists a nonzero additive function $l^*:=(\lambda^*, a^*)\in (\mathbb{R}_+)^* \times \mathcal{D}^*$, with $\lambda^*\not\equiv 0$,  such that $(a^*,\bar x)$ is a saddle point for $L(\lambda^*,\cdot,\cdot)$  on $\mathcal{D}^*\times X$.
}
\end{theorem} {\bf Proof.} By Theorem \ref{Thm:5.3+}, we have that $(0,0)$ and  $\mathcal{E}_{\bar x} $ are properly linearly separable. Following the same arguments of Theorem  \ref{Th:6.2}, we complete the proof.
\cvd

%%%%%%%%%%%%%%%%%%%%%%%%%%%%%%%%%%%%%%%%%%%%
%%%%%%%%%%%%%%%%%%%%%%%%%%%
\section{Conclusions}

We have proved several new separation theorems between convex sets in groups under the assumption that
at least one of them has nonempty pseudo relative core, whereas their quasi cores are empty, by using the  extension theorem of additive functions for groups.  This answers affirmatively the open question proposed by Li and Mastroeni in {\rm\cite{LM25}}.
We have applied our results to a class of optimization problems with cone constraints  in groups, providing Lagrangian type necessary optimality conditions. Future research may be devoted to sufficient optimality conditions for optimization problems in groups and especially in the group $\mathbb{Z}^n$,
following the line considered in \cite{LM20}.


\begin{thebibliography}{10}


%\bibitem{AF18} M. Adivar and S. C. Fang, Convex analysis and duality over discrete domains, J. Oper. Res. Soc. China,  6 (2018), pp. 189-247.

\bibitem{BP86} V. Barbu and T. Precupanu,  Convexity and Optimization in Banach Spaces, Fourth edition, Springer, Dordrecht, 2012.

\bibitem{BC17} H. H. Bauschke and P. L. Combettes,  Convex Analysis and Monotone Operator Theory in Hilbert Spaces, Second edition, Springer, Cham, 2017.

\bibitem{BG16} J. M. Borwein  and O. Giladi,  Some remarks on convex analysis in topological groups, J. Convex Anal.,  23 (2016), pp.  313-332.

\bibitem{BG18} J. M. Borwein and O. Giladi, Convex analysis in groups and semigroups: a sampler, Math. Program., Ser. B, 168 (2018), pp. 11-53.

%\bibitem{BG03} J. Borwein and R. Goebel, Notions of relative interior in Banach spaces. Optimization and related topics, J. Math. Sci. (N. Y.), 115 (2003), pp. 2542--2553.

\bibitem{BL06} J. M. Borwein and A. S. Lewis, Convex Analysis and Nonlinear Optimization. Theory and Examples, Second edition, Springer, New York, 2006.

%\bibitem{BL92} J. M. Borwein and A. S. Lewis, Partially finite convex programming. I. Quasi relative interiors and duality theory,  Math. Program., Ser. B, 57 (1992), pp. 15-48.

\bibitem{BV10} J. M. Borwein and J. D. Vanderwerff,  Convex Functions: constructions, characterizations
and counterexamples, Cambridge University Press, Cambridge, 2010.

\bibitem{Gia05} F. Giannessi, Constrained Optimization and Image Space Analysis, Springer, New York,  2005.

\bibitem{HL93} J. B. Hiriart-Urruty and  C. Lemar\'echal, Convex Analysis and Minimization Algorithms I and II, Springer-Verlag, Berlin, 1993.

\bibitem{Kau66} R. Kaufman, Interpolation of additive functionals, Studia Math., 27 (1966), pp. 269-272.

\bibitem{LM20} J. Li and G. Mastroeni,  Convex analysis in $Z^n$ and applications to integer linear programming, SIAM J. Optim., 30 (2020), pp. 2809-2840.

\bibitem{LM23} J. Li and G. Mastroeni, Optimization problems with cone constraints in groups and semigroups: an approach based on image space analysis, J. Optim. Theory Appl., 196 (2023), pp. 973-1007.

\bibitem{LM24} J. Li and G. Mastroeni, Convex analysis in groups and semigroups and applications to separation and optimization, J. Convex Anal., 32 (2025), pp. 989-1022.

\bibitem{LM25} J. Li and G. Mastroeni, Pseudo relative core of sets in groups and applications to Radon's lemma, Helly's theorem and optimization, submitted, 2026.

%\bibitem{Mat02} J. Matou\v sek, Lectures on Discrete Geometry, Springer-Verlag, New York, 2002.

%\bibitem{MR23} R. D\'iaz Mill\'an and Vera Roshchina, The intrinsic core and minimal faces of convex sets in general vector spaces, Set-Valued Var. Anal., 31 (2023), https://doi.org/10.1007/s11228-023-00671-6.

%\bibitem{MR24} R. D\'iaz Mill\'an and Vera Roshchina, Face relative interior of convex sets in topological vector spaces, SIAM J. Optim., 34 (2024), pp. 3456-3476.

%\bibitem{Mur03} K. Murota,  Discrete Convex Analysis, SIAM Monographs on Discrete Mathematics and Applications, Society for Industrial and Applied Mathematics (SIAM), Philadelphia, PA, 2003.

\bibitem{Rob96} Derek J. S. Robinson, A Course in The Theory of Groups (Second edition), Springer-Verlag, New York, 1996.

\bibitem{Roc70} R. T. Rockafellar, Convex Analysis,  Princeton University Press, Princeton, N.J., 1970.

\bibitem{rw98} R. T. Rockafellar and R. J.-B. Wets, Variational Analysis, Springer-Verlag, Berlin, 1998.

\bibitem{sw70} J. Stoer and C. Witzgall, Convexity and Optimization in Finite Dimensions  I,  Springer-Verlag, New York-Berlin, 1970.

%\bibitem{VM81} M. Volle, Regularizations in abelian complete ordered groups, J. Math. Anal. Appl., 84 (1981), pp. 418-430.

%\bibitem{VM82} M. Volle, Sur les applications regularis\'ees \'a valeurs dans un groupe  (French) [Group-valued regularized maps], Math. Operationsforsch. Statist. Ser. Optim., 13 (1982), pp. 493-500.

\bibitem{Zal02} C. Z\~alinescu, Convex Analysis in General Vector Spaces, World Scientific Publishing Co., Inc., River Edge, NJ, 2002.




\end{thebibliography}
\end{document}